\documentclass[reqno,draft]{amsart}
\usepackage{amsmath}
\allowdisplaybreaks[3]
\numberwithin{equation}{section}
\numberwithin{equation}{section}

\def\proof{\indent{\em Proof.\quad}}
\def\endproof{\hfill\hbox{$\sqcup$}\llap{\hbox{$\sqcap$}}\medskip}
\newtheorem{thm}{\indent Theorem}[section]
\newtheorem{cor}[thm]{\indent Corollary}
\newtheorem{lem}[thm]{\indent Lemma}
\newtheorem{prop}[thm]{\indent Proposition}
\newtheorem{dfn}{{\indent\bf Definition}}[section]
\newtheorem{rmk}{{\indent\bf Remark}}[section]

\newcommand{\mb}{\mbox}

\newcommand{\ol}{\overline}

\newcommand{\strl}{\stackrel}

\newcommand{\td}{\tilde}

\newcommand{\fr}{\frac}

\newcommand{\edd}{\end{document}}
\newcommand{\be}{\begin{equation}}
\newcommand{\ee}{\end{equation}}

\newcommand{\lagl}{\langle}
\newcommand{\ragl}{\rangle}
\newcommand{\lmx}{\left(\begin{matrix}}
\newcommand{\rmx}{\end{matrix}\right)}
\newcommand{\ldt}{\left|\begin{matrix}}
\newcommand{\rdt}{\end{matrix}\right|}

\newcommand{\sgn}{{\rm Sgn\,}}

\newcommand{\tr}{{\rm tr\,}}
\newcommand{\vfi}{\varphi}
\newcommand{\veps}{\varepsilon}
\newcommand{\bbr}{{\mathbb R}}

\newcommand{\bbc}{{\mathbb C}}
\newcommand{\bbh}{{\mathbb H}}
\newcommand{\bbo}{{\mathbb O}}
\newcommand{\ba}{\begin{array}}
\newcommand{\ea}{\end{array}}
\newcommand{\nnm}{\nonumber}
\newcommand{\beal}{\begin{align}}
\newcommand{\eal}{\end{align}}
\newcommand{\bea}{\begin{eqnarray}}
\newcommand{\eea}{\end{eqnarray}}

\newcommand{\id}{{\rm id}}

\newcommand{\ad}{{\rm ad}}
\newcommand{\spn}{{\rm Span\,}}

\newcommand{\ppp}[3]{\fr{\partial^2 #1}{\partial #2\partial #3}}
\newcommand{\dd}[2]{\fr{d #1}{d #2}}
\newcommand{\frkg}{{\mathfrak g}}
\newcommand{\frkh}{{\mathfrak h}}

\newcommand{\frkk}{{\mathfrak k}}
\newcommand{\frkp}{{\mathfrak p}}
\newcommand{\frkgl}{{\mathfrak g}{\mathfrak l}}
\newcommand{\frksl}{{\mathfrak s}{\mathfrak l}}

\newcommand{\cls}{{\mathcal S}}
\newcommand{\cll}{{\mathcal L}}
\newcommand{\stx}[2]{\strl{(#1)}{#2}}
\newcommand{\lalp}{\stx{\alpha}{L}\!\!_1{}}
\newcommand{\galp}{\!\!\stx{\alpha}{g}{}\!\!\!}
\newcommand{\Aalp}{\stx{\alpha}{A}{}\!\!}

\begin{document}

\title[Classification of nondegenerate equiaffine symmetric hypersurfaces]{A classification theorem of\\ nondegenerate equiaffine symmetric hypersurfaces}\thanks{Research supported by
NSFC (No. 11171091, 11371018) and partially supported by NSF of Henan Province (no. 132300410141).}%

\author[Xingxiao Li]{Xingxiao Li
}
\date{}

\begin{abstract}
Motivated by the ideas and methods used by Naitoh in the consideration of parallel totally real submanifolds in complex space forms, the author of the present paper successfully makes use of the so called Jordan triple and (restricted) structure Lie algebra associated with a given Jordan algebra to establish a one-to-one correspondence between the set of equivalence classes of connected, simply connected and nondegenerate equiaffine symmetric hypersurfaces with a given nonzero affine mean curvature and that of the equivalence classes of semi-simple real Jordan algebras. Then, via the existing classification theorem of the semi-simple real Jordan algebras with unity, a complete classification for the nondegenerate and locally equiaffine symmetric hypersurfaces with nonzero affine mean curvatures is established. As an direct application of the main theorems, we prove at the end of the paper a complete classification of nondegenerate hypersurfaces with parallel Fubini-Pick forms and nonzero affine mean curvatures.

\vskip 0.1in
{\bf Key words and expressions}\ nondegenerate hypersurfaces, equiaffine symmetric hypersurfaces, affine metric, Fubini-Pick form, Jordan algebras, pseudo-Riemannian symmetric spaces

\vskip 0.1in
{\bf 2000 AMS classification:} Primary 53A15; Secondary 53B25
\end{abstract}
\maketitle

\tableofcontents

\section{Introduction}

Let $x:M^n\to \bbr^{n+1}$ be a nondegenerate hypersurface, that is, $x$ is an immersion of which the the second fundamental form (as a Riemannian submanifolds of $\bbr^{n+1}$ endowed with the standard metric) is nondegenerate. In affine differential geometry of nondegenerate hypersurfaces, the fundamental affine invariants are the affine metric $g$ which is in general pseudo-Riemannian, the Fubini-Pick cubic form $A$ which is trilinear and symmetric, and the affine shape operator $B$. Beside these three, the affine normal $\xi:=\fr1n\Delta_gx$ and the affine mean curvature $L_1:=\tr B$ are another two important affine invariants. By lifting or lowering indices via the metric $g$, $B$ is identified with a symmetric bilinear form called the affine shape form or the affine second fundamental form; While the Fubini-Pick form $A$ can always be identified with a symmetric tensor of type $(1,2)$ (also called the difference tensor, denoted by $K$), or a ${\rm hom}(TM)$-valued $1$-form.

As all we know, hypersurfaces with particular affine invariants are the most important objects of study in the affine geometry and thus draw much attention of many authors. Examples of these are, among others (such as \cite{nom82} etc.), affine hyperspheres (with the affine shape form $B$ parallel to the affine metric $g$), affine maximal hypersurfaces (with vanishing affine mean curvature $L_1$) and those with parallel Fubini-Pick forms that are necessarily either quadratics or homogeneous affine hyperspheres (\cite{dil-vra98}, \cite{vra88}, \cite{bok-nom-sim90}).

If we only take account of the definition, affine hyperspheres seem very simple, but in fact they do form a very large class of hypersurfaces. Consequently it is a great challenge to find explicitly all the affine hyperspheres, and up to now it still remains a very hard job. In the global approach, the locally strongly convex hyperspheres have been deeply studied, see for example \cite{cal72}, \cite{ch-yau86}, \cite{amli89}, the survey \cite{lof10} and the references therein. But even assuming the global condition, it is still far away from a complete understanding of these hypersurfaces. Although this, the study of affine hyperspheres (both globally and locally) has been made a lot of great achievement by many authors. In fact, various kind of affine hyperspheres with particular properties are widely studied. For example, the proof of the Calabi's conjecture on hyperbolic hyperspheres and the rigidity theorems on the ellipsoid in the locally strongly convex case (see for example, \cite{amli89}, \cite{amli90} and \cite{amli92}); the classification of hyperspheres of constant sectional affine curvatures (\cite{vra-li-sim91}, \cite{wang93} and \cite{kri-vra99}); the study of homogeneous hyperspheres (\cite{sas80}, \cite{tsu82a}, \cite{tsu82b} for the relation with symmetric cones, \cite{dil-vra94} for the Calabi-type composition), the characterization of the Calabi composition of hyperbolic hyperspheres (\cite{hu-li-vra08}; also \cite{lix13} and \cite{lix14} in a different manner); the isotropic affine hyperspheres (\cite{bir-djo12}); and so on.

On the other hand, since the difference tensor $K$, namely, the Fubini-Pick form $A$ as a $(2,1)$ tensor measures the difference between the induced affine connection $\nabla$ and the Levi-Civita connection $\hat\nabla$ of the affine metric $g$ (in fact, $K=\nabla-\hat\nabla$), the study of the Fubini-Pick form has been of great significance and lots of interesting results has been obtained, particularly in recent several years. Among others, one classical result in this direction is the Maschke-Pick-Berwald theorem which implies that the hypersurfaces with vanishing Fubini-Pick form are necessarily quadratic affine hyperspheres (\cite{nom-sas94b},\cite{li-sim-zhao93}). Therefore, the next most natural class of hypersurfaces to consider are those with parallel Fubini-Pick form, with respect to either the induced connection $\nabla$ or the Levi-Civita connection $\hat\nabla$. Such kind of hypersurfaces have also been considered extensively and a lot of important progress in the classification issue have been made. For example, with respect to $\nabla$, the parallel of $A$ and the parallel of $K$ (the latter is not equivalent to the former in this case) are respectively considered in \cite{vra88}, \cite{nom-pin89}, \cite{nom-sas94a}, \cite{gig02}, \cite{gig03}, \cite{gig11} for $\nabla A=0$ and \cite{dil-vra98} for $\nabla K=0$. On the other hand, under the condition that $\hat\nabla A=0$ which is equivalent to that $\hat\nabla K=0$, a series of classification theorems of locally strongly convex hypersurfaces have been proved (\cite{mag-nom89}, \cite{bok-nom-sim90}, \cite{dil-vra91}, \cite{dil-vra-yap94} and \cite{hu-li-sim-vra09} for some special cases; for a complete classification, see \cite{hu-li-vra11}, and also see \cite{lix13} and \cite{lix-zhao14} in respectively two totally different manners); In \cite{amli-wang91} and \cite{wang93}, similar problems were considered respectively for the centroaffine and equiaffine hypersurfaces under the flatness assumption. As for the more general case of nondegenerate hypersurfaces, there also have been some interesting partial classification results, see for example the series of papers by Z.J. Hu et al: \cite{hu-li11}, \cite{hu-li-li-vra11a} and \cite{hu-li-li-vra11b}. In this direction, a recent development is the preprint article \cite{hil12} in which the author considers the centroaffine hypersurfaces with parallel Fubini-Pick form ($\hat\nabla A=0$) and, by an argument using the so-called $\omega$-domain and the classification of Jordan algebras, the author aims to classify all nondegenerate proper hyperspheres with $\hat\nabla A=0$, while in \cite{hil13} the same author gives a correspondence between  nondegenerate (locally) graph immersions with $\hat\nabla A=0$ and a class of real Jordan algebras admitting a nondegenerate trace form.

Since a nondegenerate hypersurface with $\hat\nabla A=0$ is not only homogeneous but also locally symmetric as a pseudo-Riemannian manifold with respect to the affine metric, in order to establish a local complete classification it is natural to introduce the concept of (equi)affine symmetric hypersurface and only needed to consider them in a global view, that is, to make a complete classification of the connected and simply connected (equi)affine symmetric ones (see \cite{lix13} and \cite{lix-zhao14}). In this present paper, motivated by the ideas and methods used by Naitoh in the consideration of parallel totally real submanifolds in complex space forms (see \cite{nai83a} and \cite{nai83b}), we make use of the so called Jordan triple (associated with a given Jordan algebra ${\mathcal J}$) to derive the structure Lie algebra $\cll$ of ${\mathcal J}$ and introduce the restricted structure Lie algebra $\frkg$. In terms of this, we successfully establish a one-to-one correspondence between the set of equivalence classes of connected, simply connected and nondegenerate equiaffine symmetric hypersurfaces with a pre-assumed affine mean curvature $L_1\neq 0$ and that of the equivalence classes of semi-simple real Jordan algebras with unity (Theorem \ref{main}). Then, via the existing classification theorem of the semi-simple real Jordan algebras with unity, a complete equiaffine classification for the locally equiaffine symmetric hypersurfaces with nonzero affine mean curvatures is established (Theorem \ref{main0}). As an direct application of the main theorems, we prove at the end of the paper a theorem that completely classifies all of those nondegenerate hypersurfaces with nonzero affine mean curvatures and parallel Fubini-Pick forms with respect to the Levi-Civita connection (Theorem \ref{thm8.3}). We should remark that, the main idea and the relevant argument we used here are completely of algebra, depending heavily on the theory of symmetric spaces (\cite{hel01}, \cite{kob-nom69}). These are totally different from those used in \cite{hil12}.

Our main theorems are stated as follows:

{\thm[The correspondence]\label{main} For a given nonzero constant $L_1$, there is a one-to-one correspondence between the set of equivalence classes of connected, simply connected and nondegenerate equiaffine symmetric hypersurfaces with affine mean curvature $L_1$ and that of the equivalence classes of semi-simple real Jordan algebras.}

{\thm[The classification]\label{main0} Let $x:M^n\to \bbr^{n+1}$ be a nondegenerate and locally equiaffine symmetric hypersurface with the affine metric $g$ and affine mean curvature $L_1\neq 0$. Then either of the following two cases should occur:
\begin{enumerate}
\item $x$ is locally affine equivalent to one of the standard imbeddings $\td x:G/K\to V$, defined by \eqref{eqn5.7}, of pseudo-Riemannian symmetric spaces $G/K$ of dimension $n$ into $(n+1)$-dimensional real simple Jordan algebras ${\mathcal J}=(V,\circ)$ where the real linear spaces $V$ are listed as follows:
\begin{enumerate}
\item $V=\bbr$, the field of real numbers, $\td x(G/K)$ contains only a point;
\item $V=\bbr^m$, the space of ordered $m$-tuple of real numbers, $n=m-1$, $m\geq 3$,
$$\td x(G/K)\subset\{(Z^tQZ)^m=C^{n+1};\ Z\in \bbr^m\},$$
where $Q$ is a symmetric $m\times m$ matrix satisfying $\det Q\neq 0$;
\item $V=M_m(\bbr)$, the space of real matrices of order $m$, $n=m^2-1$, $m\geq 3$,
$$\td x(G/K)\subset\{(\det Z)^{2m}=C^{n+1};\ Z\in M_m(\bbr)\};$$
\item $V=M_m(\bbh)$, the space of quaternion matrices of order $m$, $n=4m^2-1$, $m\geq 2$,
$$\td x(G/K)\subset\{(\det Z)^{4m}=C^{n+1};\ Z\in M_m(\bbh)\};$$
\item $V=S_m(\bbr)$, the space of real symmetric matrices of order $m$, $n=\fr12m(m+1)-1$, $m\geq 3$,
$$\td x(G/K)\subset\{(\det Z)^{m+1}=C^{n+1};\ Z\in S_m(\bbr)\};$$
\item $V=H_m(\bbc)$, the space of complex Hermitian matrices of order $m$, $n=m^2-1$, $m\geq 3$,
$$\td x(G/K)\subset\{(\det Z)^{2m}=C^{n+1};\ Z\in H_m(\bbc)\};$$
\item $V=H_m(\bbh)$, the space of quaternionic Hermitian matrices of order $m$, $n=m(2m-1)-1$, $m\geq 3$,
$$\td x(G/K)\subset\{(\det Z)^{2m-1}=C^{n+1};\ S\in H_m(\bbh)\};$$
\item $V=A_{2m}(\bbr)$, the space of real skew-symmetric matrices of order $m$, $n=m(2m-1)-1$, $m\geq 3$,
    $$\td x(G/K)\subset\{(\det Z)^{2m-1}=C^{n+1};\ Z\in A_{2m}(\bbr)\};$$
\item $V=SH_m(\bbh)$, the space of quaternionic skew-Hermitian matrices of order $m$, $n=m(2m+1)-1$, $m\geq 2$,
$$\td x(G/K)\subset\{(\det S)^{2m+1}=C^{n+1};\ S\in SH_m(\bbh)\};$$
\item $V=H_3(\bbo)$, the space of octonionic Hermitian $3\times 3$ matrices, $n=26$,
$$\td x(G/K)\subset\{\det Z=C^{\fr32};\ Z\in H_3(\bbo)\};$$
\item $V=H_3(O,\bbr)$, the space of Hermitian $3\times 3$ matrices
with entries of split octonions over $\bbr$, $n=26$,
$$\td x(G/K)\subset\{\det Z=C^{\fr32};\ Z\in H_3(O,\bbr)\};$$
\item $V=\bbc$, the field of complex numbers, $n=1$,
$$\td x(G/K)\subset\{|z|^2=C;\ z\in\bbc\};$$
\item $V=\bbc^m$, the space of ordered $m$-tuple of complex numbers, $m\geq 3$, $n=2m-1$,
$$\td x(G/K)\subset\{|Z^tZ|^{2m}=C^{n+1};\ Z\in \bbc^m\};$$
\item $V=S_m(\bbc)$, the space of complex symmetric matrices of order $m$, $n=m(m+1)-1$, $m\geq 3$,
$$\td x(G/K)\subset\{|\det Z|^{2(m+1)}=C^{n+1};\ Z\in S_m(\bbc)\};$$
\item $V=M_m(\bbc)$, the space of $m\times m$ complex matrices, $n=2m^2-1$, $m\geq 3$,
$$\td x(G/K)\subset\{|\det Z|^{4m}=C^{n+1};\ Z\in M_m(\bbc)\};$$
\item $V=A_{2m}(\bbc)$, the space of complex skew-symmetric matrices of order $m$, $n=2m(2m-1)-1$, $m\geq 3$,
    $$\td x(G/K)\subset\{|\det Z|^{2(2m-1)}=C^{n+1};\ Z\in A_{2m}(\bbc)\};$$
\item $V=H_3(O,\bbc)$, the space of Hermitian $3\times 3$ matrices with
entries of split octonions over $\bbc$, $n=53$,
$$\td x(G/K)\subset\{|\det Z|^2=C^{\fr32};\ Z\in H_3(O,\bbc)\}.$$
\end{enumerate}
where $C$ is determined by the affine mean curvature $L_1$ via
\be\label{C}C=-\sgn(L_1)\sqrt{n+1}((n+1)|L_1|)^{-\fr12(n+2)};\ee
\item $x$ is locally affine equivalent to a Calabi composition of some of the above equiaffine symmetric hypersurfaces including the $0$-dimensional ones.
\end{enumerate}}

The present paper is organized as follows:

In the next section (Section \ref{sec2}), we briefly review some relevant material on equiaffine differential geometry of nondegenerate hypersurfaces and the basics on the Jordan algebra, including a classification theorem for the semi-simple Jordan algebra.

In Section \ref{sec3}, for a given connected, simply connected and nondegenerate equiaffine symmetric hypersurface $x:M^n\to\bbr^{n+1}$ with affine metric $g$, Fubini-Pick form $A$ and affine mean curvature $L_1$, we first explicitly construct a symmetric pair of Lie algebras $(G,K)$ such that $M^n=G/K$, which is needed in the final classification argument. Then we fix one point $o\in M^n$ and use the invariance property of the Fubini-Pick form $A$ to define a unital Jordan algebra ${\mathcal J}=(V,\circ)$ with $V=x_{*o}(T_oM^n)\oplus\bbr\cdot x(o)\equiv\bbr^{n+1}$, where the Jordan algebra product $\circ$ is essentially defined by
$$
x_{*o}(X)\circ x_{*o}(Y)=x_{*o}(A(X,Y))-L_1g(X,Y)e,\quad \forall X,Y\in T_oM^n,
$$
where $e=C^{-1}x(o)$ is taken as the unity of ${\mathcal J}$. We show that ${\mathcal J}$ is semi-simple if and only if $L_1\neq 0$, and different choices of the given point $o$ define isomorphic Jordan algebras.

In Section \ref{sec4}, by making use of the Jordan triple associated with a given semi-simple Jordan algebra ${\mathcal J}=(V,\circ)$ which is necessarily unital, we explicitly present the so called structure Lie algebra $\cll$ and derive an effective symmetric pair $(\frkg,\frkk)$ of Lie algebras (Proposition \ref{prop4.1}), where $\frkg$ is the restricted structure Lie algebra of ${\mathcal J}$.

In section \ref{sec5}, for any given semi-simple Jordan algebra ${\mathcal J}=(V,\circ)$, we first find the connected and simply connected symmetric space $G/K$ that corresponds to the symmetric pair $(\frkg,\frkk)$ of Lie algebras derived in Section \ref{sec4}; Then we define for each $L_1\neq 0$ an equivariant imbedding $x:G/K\to V$ which is proven a connected, simply connected and nondegenerate equiaffine symmetric hypersurface with parallel Fubini-Pick form ($\hat\nabla A=0$) and affine mean curvature $L_1$ (Proposition \ref{prop5.2}). We also find the defining equation of the image $x(G/K)$ (Proposition \ref{prop5.3}).

In section \ref{sec6}, we give the proof of Theorem \ref{main}.

In section \ref{sec7}, we introduce and discuss in detail one Calabi-type composition of nondegenerate equiaffine symmetric hypersurfaces. These are needed in the argument for the main classification theorem.

In the last section (Section \ref{sec8}) we prove the main classification theorem (Theorem \ref{main0}) by using Theorem \ref{main} and the discussion in Section \ref{sec7}. As a direct application of the main theorems, we arrive at a complete classification of all nondegenerate hypersurfaces with parallel Fubini-Pick form and nonzero affine mean curvature (see Theorem \ref{thm8.3}).

{\sc Acknowledgement} The first author is grateful to Professor A-M Li for his encouragement and important suggestions during the preparation of this article. He also thanks Professor T. Sasaki and Professor Z. J. Hu for providing him valuable related references some of which are listed in the end of this paper. Moreover, the very first result of the author in the approach of affine symmetric hypersurfaces was announced at the international conference on differential geometry (Bedlewo, Poland, June, 2012). The author would like to express his hearty thanks to the conference organizers, in particular, Professor B. Opozda and Professor U. Simon for their kind invitation and hospitality.

\section{Preliminaries}\label{sec2}

\subsection{The equiaffine geometry of hypersurfaces}

In this subsection, we brief some basic facts in the equiaffine geometry of nondegenerate hypersurfaces. For details the readers are referred to, say, \cite{nom82}, \cite{li-sim-zhao93} and \cite{nom-sas94b}.

Let $x:M^n\to\bbr^{n+1}$ be a nondegenerate hypersurface. For any vector field $\eta\in\bbr^{n+1}$ transversal to the tangent space $x_*(TM^n)$ of $x$, we have the following direct decomposition of vector spaces
$$
x^*T\bbr^{n+1}=x_*(TM^n)+\bbr\cdot\eta.
$$
This decomposition and the canonical differentiation $\bar D^0$ on $\bbr^{n+1}$ define a nondegenerate bilinear form $h\in\bigodot^2T^*M^n$ and a connection $D^\eta$ on $TM^n$ as follows:
\be\label{dfn h}
\bar D^0_Xx_*(Y)=x_*(D^\eta_XY)+h(X,Y)\eta,\quad\forall X,Y\in TM^n.
\ee
\eqref{dfn h} is sometimes called the affine Gauss formula of the hypersurface $x$.

In what follows we make the following convention for the range of indices:
$$1\leq i,j,k,l\leq n.$$

Let $\{e_i,e_{n+1}\}$ be a local unimodular frame field along $x$ with $e_i\in TM$ and with $\{\omega^i,\omega^{n+1}\}$ its dual coframe.  Then $\eta:=e_{n+1}$ is transversal to the tangent space $x_*(TM^n)$. Write $h=\sum h_{ij}\omega^i\omega^j$ with $h_{ij}=h(e_i,e_j)$ and $H=|\det(h_{ij})|$. Then (see \cite{li-sim-zhao93}) the locally defined nondegenerate metric
\be\label{affine metric}g:=H^{-\fr1{n+2}}h\ee
is independent of the choice of the unimodular frame field $\{e_i,e_{n+1}\}$ and thus is in fact a globally well-defined metric on $M^n$ which is called the affine (or Berwald-Blaschke) metric. By taking $x$ as an $\bbr^{n+1}$-valued smooth function on $M^n$, we call the vector function
\be\label{affine normal}\xi:=\fr1n\tr_g(x)\ee
the affine normal vector.

{\rmk\rm When $x$ is locally strongly convex, the bilinear form $h$ is either positive definite or negative definite. Thus we can choose suitable transversal vectors $e_{n+1}$ (or equivalently, choose a suitable orientation of $M^n$) such that the affine metric $g$ defined in \eqref{affine metric} is positive definite.}

If the transversal vector field $\eta$ in \eqref{dfn h} is chosen to be parallel to the affine normal $\xi$, then the induced connection $\nabla:=D^\eta$ is independent of the choice of $\eta$ and is called the induced affine connection of $x$.

If $\hat\nabla$ is the Levi-Civita connection of the affine metric $g$, then
the difference tensor $K$, namely, the Fubini-Pick form (as a symmetric $(1,2)$ tensor) is defined by
\be\label{f-p}
K(X,Y)\equiv A(X,Y)=\nabla_XY-\hat\nabla_XY,\quad \forall\, X,Y\in TM^n,
\ee
which is identified via the affine metric $g$ with a symmetric cubic form $A(X,Y,Z)=g(A(X,Y),Z)$. This cubic form $A$ is also called the Fubini-Pick form (or simply the cubic form) of $x$. For each $X\in TM^n$, we have a linear map $A_X:TM^n\to TM^n$ given by
\be\label{A_X}A_X(Y):=A(X,Y),\quad \forall Y\in TM^n.\ee

The affine shape operator $B:TM^n\to TM^n$ of $x$ is defined by $X(\xi)=-x_*(B(X))$ for $X\in TM$. Via the affine metric $g$, $B$ corresponds to a symmetric bilinear form $B\in\bigodot^2T^*M^n$ by $$B(X,Y)=g(B(X),Y),\quad X,Y\in TM^n,$$
which is called the equiaffine Weingarten form (\cite{li-sim-zhao93}) or the (equi)affine second fundamental form of $x$.

It is known that the assumption that $\eta\parallel\xi$ is equivalent to the apolarity of $A$, that is
\be\label{apol}
\tr A_X=0,\text{\ or equivalently\ }\sum_{i,j}g^{ij}A(X,e_i,e_j)=0,\quad \forall X\in TM^n,
\ee
where the matrix $(g^{ij})=(g_{ij})^{-1}$ and $g_{ij}=g(e_i,e_j)$.

In this sense, the affine Gauss equation can be written as follows:
\be\label{gaus}
R(X,Y)Z=\fr12(g(Y,Z)B(X)+B(Y,Z)X-g(X,Z)B(Y)-B(X,Z)Y)-[A_X,A_Y](Z).
\ee

Furthermore, it is not hard to show that the Fubini-Pick form $A$, the affine second fundamental form $B$ and the affine mean curvature $L_1$ meet the following equation:
\be\label{abl1}
B=L_1g+\fr2n{\rm div}(A), \text{\ where\ } ({\rm div}A)(X,Y)=\omega^i(\hat\nabla_{e_i}A)(X,Y),\quad \forall X,Y\in TM^n.
\ee

Each of the eigenvalues $B_1,\cdots,B_n$ of the affine shape operator $B$ is called the affine principal curvature of $x$. Moreover, the the affine mean curvature of $x$ is defined by
\be\label{afme}
L_1:=\fr1n\tr B=\fr1n\sum_iB_i.
\ee
A hypersurface $x$ is called an proper (resp. improper) affine hypersphere, if all of its affine principal curvatures are equal to a nonzero (resp. zero) constant. In such case we have
\be\label{afsp}
B=L_1g,\ \text{or equivalently,\ } B(X)=L_1X \quad\mb{for all\ }X\in TM^n.
\ee
It follows that the affine Gauss equation \eqref{gaus} of an affine hypersphere assumes the following form:
\be\label{gaus_af sph}
R(X,Y)Z=L_1(g(Y,Z)X-g(X,Z)Y)-[A_X,A_Y](Z),
\ee

Furthermore, all the affine normal lines of a proper affine hypersphere $x:M^n\to\bbr^{n+1}$ pass through a fix point $O$ called the affine center of $x$; while the affine normals of an improper affine hypersphere are parallel to a fixed direction.

{\rmk\rm When $x$ is locally strongly convex, a suitable orientation is chosen to ensure that the affine metric $g$ is positive definite. It follows that the sign of the affine mean curvature $L_1$ is an equiaffine invariant. Thus in this case the affine hyperspheres are divided into three classes: elliptic affine hyperspheres ($L_1>0$), parabolic affine hyperspheres ($L_1=0$) and hyperbolic affine hyperspheres ($L_1<0$)}.

The following affine uniqueness theorem is needed later:

{\thm\label{af uniq} $($\cite{li-sim-zhao93}$)$ Let $x:M^n\to \bbr^{n+1}$,
$\bar x:\bar M^n\to \bbr^{n+1}$ be two nondegenerate hypersurfaces with respectively the affine metrics $g$, $\bar g$ and the Fubini-Pick forms $A$, $\bar A$. Let $\vfi:(M^n,g)\to (\bar M^n,\bar g)$ be an isometry between pseudo-Riemannian manifolds. Then $\vfi^*\bar A=A$ if and only if there exists a unimodular affine transformation $\Phi:\bbr^{n+1}\to \bbr^{n+1}$ such that $\bar x\circ\vfi=\Phi\circ x$, or equivalently, $\bar x=\Phi\circ x\circ\vfi^{-1}$.}

\proof For those locally strongly convex hypersurfaces, the theorem is proved in \cite{li-sim-zhao93}. A careful examination shows that the same argument in \cite{li-sim-zhao93} also applies for the general nondegenerate case.\endproof

{\rmk\rm The proof of the sufficient part of Theorem \ref{af uniq} can also be found in \cite{lix13} and \cite{lix14}.}

Motivated by Theorem \ref{af uniq}, we introduce the following concept of affine equivalence relation between nondegenerate hypersurfaces:

{\dfn\rm Let $x:M^n\to \bbr^{n+1}$ be a nondegenerate hypersurface with the affine metric $g$. A hypersurface $\bar x:M^n\to \bbr^{n+1}$ is called affine equivalent to $x$ if there exists a unimodular transformation $\Phi:\bbr^{n+1}\to \bbr^{n+1}$ and an isometry $\vfi$ of $(M^n,g)$ such that $\bar x=\Phi\circ x\circ\vfi^{-1}$.}

To end this section, we would like to recall the following definition:

{\dfn\label{dfn afsym} \rm{(\cite{lix13})} A nondegenerate hypersurface $x:M^n\to \bbr^{n+1}$ is called equiaffine symmetric (resp. locally equiaffine symmetric) if

$(1)$ the pseudo-Riemannian manifold $(M^n,g)$ is symmetric (resp. locally symmetric);

$(2)$ the (equiaffine) Fubini-Pick form $A$ is invariant (resp. locally invariant) on $M^n$.}

\subsection{The real Jordan algebras}

Jordan algebras were first introduced by Pascual Jordan in 1933 to formalize the notion of an algebra of observables in quantum mechanics (see for example \cite{jor-neu-wig34}). They were originally called "r-number systems", but were renamed "Jordan algebras" by Albert (1946), who began the systematic study of general Jordan algebras. For later use, we make a brief review at this moment here of the Jordan algebras over real numbers and complex numbers, including particularly the related classification theorems. The main references in this subsection are \cite{jac68}, \cite{koc99}, \cite{mcc04} and \cite{hil12}.

{\dfn\rm Let $(V,\circ)$ be a finite dimensional algebra over the field ${\bbr}$ of real numbers (resp. the field ${\mathbb C}$ of complex numbers). Then $(V,\circ)$ is called a real (resp. complex) Jordan algebra if the following two conditions are satisfied:

(JA1: Commutativity) For all $u,v\in V$, $u\circ v=v\circ u$;

(JA2: Jordan Identity) For all $u,v\in V$, $u\circ(u^2\circ v)=u^2\circ(u\circ v)$

\noindent
where for any integer $k\geq 2$, $u^k=u\circ u^{k-1}$. If it is the case, the product $\circ$ is called the Jordan product.}

In what follows we always suppose that ${\mathcal J}:=(V,\circ)$ is a real or complex Jordan algebra. For each $v\in V$, there correspond two linear maps $T_v,P_v:V\to V$ by $T_v(u)=v\circ u$, $u\in V$ and $P_v:=2T_v^2-T_{v^2}$. Thus it is natural to define the determinant and the trace of any $v\in V$ by
\be\label{eqn2.17}\det v=\det P_v,\quad\tr v=\tr T_v\ee respectively. Define
\be\label{eqn2.18}
\lagl u,v\ragl=\tr T_{u\circ v},\quad \forall u,v\in V.
\ee
Then $\lagl \cdot,\cdot\ragl$ is a symmetric bilinear form on $V$. In some cases we shall write $\lagl\cdot,\cdot\ragl\!_{\mathcal J}$ for $\lagl\cdot,\cdot\ragl$ if necessary.

{\dfn\label{dfn24}\rm Let ${\mathcal J}=(V,\circ)$ be a real or complex Jordan algebra.

(1) ${\mathcal J}$ is called nondegenerate if the linear map $T_v:V\to V$ is nondegenerate for each $v\in V$;

(2) ${\mathcal J}$ is said to be semi-simple if the symmetric form $\lagl \cdot,\cdot\ragl$ is nondegenerate;

(3) ${\mathcal J}$ is called the direct sum of non-trivial Jordan subalgebras if there are non-trivial Jordan subalgebras ${\mathcal J}_i=(V_i,\circ)$, $i=1,\cdots,r$ ($r\geq 1$), such that
$$
V=V_1\oplus\cdots\oplus V_r,\quad V_i\circ V_j=\{0\}\text{\ for\ }i\neq j;
$$
In this case, all ${\mathcal J}_i$'s are necessarily ideas of ${\mathcal J}$.

(4) ${\mathcal J}$ is called simple if it is semi-simple and can not be decomposed into a direct sum of its non-trivial Jordan subalgebras, that is, ${\mathcal J}$ has no non-trivial Jordan ideas;

(5) An element $v\in V$ is called nondegenerate if the associated linear map $P_v$ is nondegenerate, or equivalently, $\det P_v\neq 0$.

(6) ${\mathcal J}$ is said to be unital if there is an element $e\in V$, called the unity, such that $v\circ e=v$ for all $v\in V$; In this case, an element $v\in V$ is said to be invertible if there exists some $v^{-1}\in V$, called the inverse of $v$, such that $v\circ v^{-1}=e$ and $[T_v,T_{v^{-1}}]=0$;

(7) ${\mathcal J}$ is said to be central-simple if its center
$$
{\mathfrak z}({\mathcal J}):=\{v\in V;\ [T_v,T_u]=0\text{\ for all\ }u\in V\}
$$
is of dimension $1$.}

{\prop\label{prop2.2} (\cite{koc99}) Let J be a semi-simple real Jordan algebra. Then there exist simple subalgebras ${\mathcal J}_1,\cdots, {\mathcal J}_r$ of ${\mathcal J}$ such that ${\mathcal J}={\mathcal J}_1\oplus\cdots\oplus{\mathcal J}_r$; Furthermore,
any such two decompositions of ${\mathcal J}$ are equal up to permutation of the summands.}

{\prop\label{prop2.3} (\cite{koc99}) Every semi-simple Jordan algebra is unital.}

{\prop\label{prop2.5} $($\cite{koc99}$)$ The linear maps $T_v$ and $P_v$, $v\in V$, are self-adjoint with respect to the symmetric bilinear form $\lagl\cdot,\cdot\ragl$.}

{\prop\label{prop2.6} $($\cite{koc99}$)$ Let ${\mathcal J}=(V,\circ)$ be unital. Then the following conclusions hold:

(1) An element $v\in V$ is invertible if and only if it is nondegenerate;

(2) For an invertible element $v\in V$,
\begin{align}
&P_{v^{-1}}=P_v^{-1},\text{\ implying that\ }\det P_{v^{-1}}=(\det P_v)^{-1},\\
&v^{-1}=P_v^{-1}v,\text{\ implying that\ }T_{v^{-1}}=T_v P_v^{-1}=P_v^{-1}T_v;
\end{align}

(3) $P_{P_uv}=P_u\circ P_v\circ P_u$ for all $u,v\in V$, where $\circ$ denotes the composition of transformations.}

{\prop\label{prop2.7} $($\cite{jac68}$)$ A simple real Jordan algebra ${\mathcal J}:=(V,\circ)$ is either central-simple or it is (isomorphic to) a central-simple complex Jordan algebra viewed as a real one.}

For a give unital Jordan algebra ${\mathcal J}=(V,\circ)$ and a fixed element $\Gamma\in V$. Define a
new product $\circ_\Gamma$ on the linear space $V$ by
\be\label{isot-prod}
u\circ\!_\Gamma v=u\circ(v\circ\Gamma)+v\circ(u\circ\Gamma)-(u\circ v)\circ \Gamma.
\ee
Let $T^\Gamma_u,P^\Gamma_u\in\frkgl(V)$ for $u\in V$ be defined similar to those given
previously for ${\mathcal J}$, that is,
$$
T^\Gamma_u(v)=u\circ\!_\Gamma v,\ \forall v\in V;\quad
P^\Gamma_u=2(T^\Gamma_u)^2-T^\Gamma_{u\circ_\Gamma u}.
$$

Clearly, $T^\Gamma_u=[T_u,T_\Gamma]+T_{u\circ \Gamma}$.

{\thm\label{isot-thm} $($\cite{koc99}$)$ The new product $\circ_\Gamma$ is a Jordan product on $V$ which makes ${\mathcal J}_\Gamma:=(V,\circ_\Gamma)$ a new Jordan algebra. Furthermore,

(1) ${\mathcal J}_\Gamma$ is unital if and only if $\Gamma$ is invertible. In this case, the unity
of ${\mathcal J}_\Gamma$ is $e^\Gamma=\Gamma^{-1}$;

(2) If ${\mathcal J}_\Gamma$ is semi-simple, then ${\mathcal J}_\Gamma$ is semi-simple if and only if $\Gamma$ is invertible.}

{\dfn\rm In case that $\Gamma\in V$ is invertible, the Jordan algebra ${\mathcal J}_\Gamma$ is called
the $\Gamma$-isotope of ${\mathcal J}$ or a mutation of ${\mathcal J}$ with respect to $\Gamma$.}

\rmk\label{twisted}\rm There are some unital Jordan algebras ${\mathcal J}=(V,\circ)$ with $V$ consisting
of certain square matrices of entries of either real, complex numbers, or quaternions, or
octonions. In such cases, we can particularly choose a diagonal matrix $\Gamma\in V$ whose
diagonal elements are $\pm 1$. Then the corresponding $\Gamma$-isotope ${\mathcal J}_\Gamma$ is
referred to as the {\it twisted algebra} of ${\mathcal J}$ respect to $\Gamma$.

{\thm\label{thm2.8} Let ${\mathcal J}=(V,\circ)$ be a simple real Jordan algebra. Then

\begin{enumerate}
\item ${\mathcal J}$ is central-simple and linear isomorphic to one of the following real Jordan algebras(\cite{jac68}):
\begin{enumerate}
\item $\bbr$, the field of real numers;
\item $Jord_m(Q_\bbr)$, a real quadratic factor of dimension $m\geq 3$;
\item $M_m(\bbr)$, the Jordan algebra of real square matrices of order $m\geq 3$;
\item $M_m(\bbh)$, the Jordan algebra of quaternion square matrices of order $m\geq 2$;
\item $S_m(\bbr,\Gamma)$, the twisted algebra of real symmetric matrices of order $m\geq 3$;
\item $H_m(\bbc,\Gamma)$, the twisted algebra of complex Hermitian matrices of order $m\geq 3$;
\item $H_m(\bbh,\Gamma)$, the twisted algebra of quaternionic Hermitian matrices of order $m\geq 3$;
\item $H_m(Q,\bbr)$, the real Jordan algebra of Hermitian matrices of order $m\geq 3$ with entries of split quaternion over $\bbr$;
\item $SH_m(\bbh)$, the real Jordan algebra of skew Hermitian quaternion matrices of order $m\geq 2$;
\item $H_3(\bbo,\Gamma)$, the twisted algebra of octonionic Hermitian matrices of order $3$;
\item $H_3(O,\bbr)$, the real Jordan algebra of Hermitian matrices of order $3$ with entries of split octonion over $\bbr$.
\end{enumerate}
or
\item ${\mathcal J}$ is not central-simple and linear isomorphic to the following complex Jordan algebras viewed as real ones (\cite{mcc04}):
\begin{enumerate}
\item $\bbc$, the field of complex numbers;
\item $Jord_m(I)$, the complex quadratic factor of dimension $m\geq 3$;
\item $S_m(\bbc)$, the Jordan algebra of complex symmetric matrices of order $m\geq 3$;
\item $M_m(\bbc)$, the Jordan algebra of complex square matrices of order $m\geq 3$;
\item $H_m(Q,\bbc)$, the complex Jordan algebra of Hermitian matrices of order $m\geq 3$ with entries of split quaternion over $\bbc$;
\item $H_3(O,\bbc)$, the complex Jordan algebra of Hermitian matrices of order $3$ with entries of split octonion over $\bbc$.
\end{enumerate}
\end{enumerate}}

{\rmk\rm Theorem \ref{thm2.8} can be also found in \cite{hil12}.}

\section{The real Jordan algebras from nondegenerate equiaffine symmetric hypersurfaces}\label{sec3}

Let $x:M^n\to\bbr^{n+1}$ be a connected, simply connected and nondegenerate equiaffine symmetric hypersurface with affine metric $g$, Fubini-Pick form $A$ and affine mean curvature $L_1$. Then by Definition \ref{dfn afsym}, $(M^n,g)$ is a connected, simply connected pseudo-Riemannian symmetric space which can always be identified with some homogeneous space $G/K$ with $G$ connected and simply connected and with $K\subset G$ connected and closed. In particular, the curvature tensor $R$ is parallel with respect to the Levi-Civita connection of $g$. For later use we need to provide a rather explicit construction of the symmetric pair $(\frkg,\frkk)$ of Lie algebras for the symmetric pair $(G,K)$ of Lie groups.

Note that, in general, a symmetric pair $(\frkg,\frkk)$ is a finite dimensional real Lie algebra $\frkg$ together with a Lie subalgebra $\frkk\subset\frkg$ such that $\frkg$ can be decomposed as
$\frkg=\frkk\oplus\frkp$ where $\frkp\subset \frkg$ is a linear subspace, such that
$$
[\frkk,\frkp]\subset\frkp,\quad [\frkp,\frkp]\subset \frkk.
$$
A pair $(\frkg,\frkk)$ of Lie algebras is called effective if $\frkk$ intersects the center ${\mathfrak z}(\frkg)$ of $\frkg$ trivially, that is, $\frkk\cap{\mathfrak z}(\frkg)=\{0\}$. It is known that a symmetric pair $(\frkg,\frkk)$ is equivalent to a symmetric Lie algebra which is a pair $(\frkg,\theta)$ where $\theta$ is an involutive automorphism on $\frkg$, that is, $\theta:\frkg\to\frkg$ is an automorphism satisfying $\theta^2=\id_\frkg$, $\theta\neq \id_{\mathfrak g}$. They are related by
\be\label{eqn4.10a}
{\mathfrak k}=\{X\in {\mathfrak g};\ \theta(X)=X\},\quad
{\mathfrak p}=\{X\in {\mathfrak g};\ \theta(X)=-X\}.
\ee

To construct the pair $(\frkg,\frkk)$, we fix one point $o\in M^n$. Set $V_0=T_oM^n$ and $V=V_0\oplus \bbr\cdot e\equiv\bbr^{n+1}$ where $e:=C^{-1}x(o)$ with $C$ defined by \eqref{C}.

For each $u\in V$, there associates a linear map $T_u:V\to V$ defined by
\begin{align}
(1)& \text{\ for\ }u=\lambda e,\  T_u(v)=\lambda v,\ \forall\lambda\in\bbr \text{\ and\ }\forall v\in V;\label{eqn3.1}\\
(2)& \text{\ for\ }u=X\in V_0,\ T_X(\lambda e)=\lambda X,\ T_X(Y)=A(X,Y)-L_1g(X,Y)e,\quad\forall\lambda\in\bbr,\ \forall Y\in V_0.\label{eqn3.2}
\end{align}

Denote by $\frkgl(V)$ and $\frksl(V)$, respectively, the general linear Lie algebra and the special linear Lie algebra on the linear space $V$. Then there is a linear map $T:V\to \frkgl(V)$ which sends each $u\in V$ to $T_u\in \frkgl(V)$. Clearly, $T_u(v)=T_v(u)$ for all $u,v\in V$ and the kernel of $T$ is trivial so we can identify $V$ with its image $T(V)$ in $\frkgl(V)$; In particular, we can identify $V_0$ with the subspace $T(V_0)\subset \frkgl(V)$. For any $X,Y\in V_0$, extend the curvature tensor $R(X,Y)$ at $o$ to be an element of $\frkgl(V)$ by $R(X,Y)(e)=0$. In this sense, define
\be\label{eqn3.3}
\frkk={\rm gen}\{R(X,Y);\ X,Y\in V_0\}\subset\frkgl(V),
\ee
where ``${\rm gen}\{\cdots\}$'' denotes the Lie algebra generated by the subset ``$\{\cdots\}$''. Then we have

{\lem\label{lem3.1} $\frkk$ is a Lie subalgebra of $\frkgl(V)$. Furthermore, $\frkk$ annihilates both the affine metric $g$ and its curvature tensor $R$ at $o$.}

\proof By definition, $\frkk$ is clearly a linear subspace of $\frkgl(V)$. Since $(M^n,g)$ is a pseudo-Riemannian symmetric space, the actions of $R(X,Y)$ on the affine metric $g$ and the curvature tensor $R$ vanish for all $X,Y\in V_0$, that is
\begin{align}
(R(X,Y)\cdot g)(Z,W)\equiv& -g(R(X,Y)Z,W)-g(Z,R(X,Y)W)=0,\\
(R(X,Y)\cdot R)(Z,W)\equiv& [R(X,Y],R(Z,W)]-R(R(X,Y)Z,W)-R(Z,R(X,Y)W)=0,\label{eqn3.3-0}\\
&\quad\forall Z,W\in V_0.
\end{align}
These two equalities exactly mean that $\frkk$ annihilates both the affine metric $g$ and its curvature tensor $R$ at $o$. Furthermore \eqref{eqn3.3-0} also indicates that the bracket for linear transformations is closed in $\frkk$ and hence $\frkk$ is a Lie subalgebra of $\frkgl(V)$.
\endproof

{\rmk\rm The Lie algebra $\frkk$ is in fact the holonomy Lie algebra due to the Ambrose-Singer theorem (\cite{kob-nom63} or directly \cite{amb-sin53}).}

Next we define
\be\label{eqn3.3-1}\frkp=T(V_0)\equiv V_0,\quad \frkg=\frkk\oplus\frkp,\ee
and the bracket product on $\frkg$ by
\be\label{eqn3.4}
[T_X,T_Y]=-R(X,Y),\quad \quad [\Phi,T_X]=T_{\Phi(X)},\quad [\Phi,\Psi]=\Phi\circ\Psi-\Psi\circ\Phi,
\ee
for all $X,Y\in V_0$ and $\Phi,\Psi\in \frkk$, where ``$\circ$'' denotes the composition of transformations.

{\lem With respect to the bracket $[\cdot,\cdot]$ defined in \eqref{eqn3.4}, $\frkg$ is a Lie algebra, and $\frkk$ is a Lie subalgebra of $\frkg$.}

\proof To prove that $\frkg$ is a Lie algebra, it suffices to check the Jacobi identity
\be\label{jac-id}
[Z_1,[Z_2,Z_3]]+[Z_2,[Z_3,Z_1]]+[Z_3,[Z_1,Z_2]]=0,\quad\forall Z_1,Z_2,Z_3\in\frkg,
\ee
since the bracket defined by \eqref{eqn3.4} is clearly bilinear and skew-symmetric.

If all these $Z_i$'s belong to $\frkk$, then \eqref{jac-id} is nothing but the Jacobi identity for $\frkk$; If $Z_1,Z_2\in \frkk$ and $Z_3\in\frkp$, then \eqref{jac-id} is immediate from the definition of the bracket; If $Z_1,Z_2\in\frkp$ and $Z_3\in\frkk$, then \eqref{jac-id} reduces to the identity \eqref{eqn3.3-0}; Finally if $Z_1,Z_2,Z_3\in \frkp$, then \eqref{jac-id} is given by the Bianchi identity:
\begin{align*}
[T_{Z_1},[T_{Z_2},T_{Z_3}]]+&[T_{Z_2},[T_{Z_3},T_{Z_1}]]+[T_{Z_3},[T_{Z_1},T_{Z_2}]]\\ =&-[Z_1,R(Z_2,Z_3)]-[Z_2,R(Z_3,Z_1)]-[Z_3,R(Z_1,Z_2)]\\
=&[R(Z_2,Z_3)Z_1+R(Z_3,Z_1)Z_2+R(Z_1,Z_2)Z_3=0.
\end{align*}
\endproof

Let $\theta:\frkg\to \frkg$ be defined by $\theta(\Phi+T_X)=\Phi-T_X$ for all $\Phi\in\frkk$, $X\in V_0$. Then $\theta$ is clearly an involutive automorphism of the Lie algebra $\frkg$, that is, $\theta^2=\id_\frkg$ and $\theta\neq \id_\frkg$.

{\lem\label{lem3.3} The $(\frkg,\frkk)$ is an effective symmetric pair of Lie algebras, or equivalently, $(\frkg,\theta)$ is an effective symmetric Lie algebra.}

\proof Clearly, $(\frkg,\frkk)$ is a symmetric pair of Lie subalgebras by \eqref{eqn3.3-1} and \eqref{eqn3.4}. Suppose that $\frkh$ is an ideal of $\frkg$ which is contained in $\frkk$. If $\Phi\in\frkh$, $T_X\in\frkp$ ($X\in V_0$), then
$$
T_{\Phi(X)}=[\Phi,T_X]\in \frkp\cap\frkh=\{0\},
$$
that is, $T_{\Phi(X)}=0$ implying $\Phi(X)=0$ for all $X\in V_0$ since the Jordan algebra ${\mathcal J}$ is nondegenerate. But $\Phi(e)=0$ by the definition of extension, we have $\Phi(u)=0$ for all $u\in V$ which implies that, as an element of $\frkgl(V)$, $\Phi=0$. It follows that $\frkh=\{0\}$.
\endproof

{\prop\label{prop3.3} There exists a connected and simply connected pseudo-Riemannian symmetric space $G/K$ corresponding to the symmetric pair $(\frkg,\frkk)$ of Lie algebras. Moreover, $M^n$ together with the affine metric $g$ is isometric to $G/K$.}

\proof Let $G$ be the connected and simply connected Lie group with Lie algebra $\frkg$. Then by the simply connectedness of $G$, the involution $\theta$ lifts to an involutive automorphism $\sigma:G\to G$ (cf. \cite{che46}, p.113). Let $K$ be the unit component of the subgroup $K_\sigma$ of $G$ defined by
$$K_\sigma=\{\vfi\in G;\ \sigma(\vfi)=\vfi\}.$$ Then $K$ is closed (since $K_\sigma$ is) and with $\frkk$ as its Lie algebra, and $G/K$ is a connected and simply connected symmetric space of which the symmetry $s_o$ at the origin $o$ is given by $s_o(gK)=\sigma(g)K$. $G/K$ is almost effective since $(\frkg,\frkk,\theta)$ is effective. Generally $G/K$ may not be effective, but without loss of generality we can assume it is.

On the other hand, by Lemma \ref{lem3.1}, $\frkk$ annihilates both the affine metric $g$ and its curvature tensor $R$ at $o$. This means that $g$ and $R$ at $o$ are $ad_{\frkg}(\frkk)$-invariant and thus are $Ad(K)$-invariant. Therefore, there exists a $G$-invariant metric $\td g$ with $G$-invariant curvature tensor $\td R$ such that at the origin $o$ they coincide with $g$ and $R$ respectively. Note that both tangent spaces $T_oM^n$ and $T_oG/K$ are identified with $\frkp\equiv V_0$. Then the second conclusion of the proposition follows easily from the simply connectedness of $M^n$, $G/K$ and the following two lemmas one of which can be found directly in \cite{hel01} and the other is a slight modification (including the proof) of a lemma in \cite{hel01}. \endproof

{\lem (\cite{hel01}, Page 200, Lemma 1.2) Let $M$ and $M'$ be two manifolds with affine connections $\nabla$ and $\nabla'$, respectively, $p\in M$, $p'\in M'$, and $A:T_pM\to T_{p'}M'$ be a linear isomorphism. Assume that the corresponding torsion tensors and curvature tensors $T,T',R,R'$ are all parallel and $A\cdot T=T'$, $A\cdot R=R'$. Then there exist some open subsets $U_p\subset M$, $U'_{p'}\subset M'$ and affine diffeomorphism $\vfi:U_p\to U'_{p'}$ such that $\vfi(p)=p'$ and $\vfi_{*p}=A$.}

{\lem (cf. \cite{hel01}, Page 201, Lemma 1.4) Let $\vfi:(M,g)\to (M',g')$ be an affine diffeomorphism between two pseudo-Riemannian manifolds. If there exists some point $p\in M$ such that $\vfi_p:T_pM\to T_{\vfi(p)}M'$ is a linear isometry, then $\vfi$ is an isometry of $(M,g)$ onto $(M',g')$.}

{\rmk\label{rmk3.1}\rm Without loss of generality, in what follows we simply identify $M^n$ with $G/K$. It then follows from the parallel of $g$ and Definition \ref{dfn afsym}
that $G$ keeps invariant both the affine metric $g$ and the Fubini-Pick form $A$. This fact together with the affine uniqueness theorem (Theorem \ref{af uniq}) shows that $G$ can be identified with some subgroup of the unimodular transformation  group on $\bbr^{n+1}$.}

Now, as the main step, we introduce a Jordan algebra associated with the hypersurface $x:M^n\to\bbr^{n+1}$. To this end, the following conclusion is needed:

{\prop\label{prop3.6} Let $A$ be a covariant tensor field on a pseudo-Riemannian symmetric space $G/K$. Then $A$ is $G$-invariant if and only if it is parallel with respect to the Levi-Civita connection.}

\proof To simplify matters, we can assume without loss of generality that $A$ is a covariant tensor field of order two. Let $(\mathfrak{g},\mathfrak{k})$ be the symmetric pair corresponding to the symmetric space $G/K$ and ${\mathfrak g}=\mathfrak{k}+\frkp $ be the canonical decomposition of the pair $(\mathfrak{g},\mathfrak{k})$. Then the vector space $\frkp $ is identified with $T_oM$. Let $o=eK\in M^n$ be the base point with $e$ the identity of $G$.
Note that, for all $X,Y,Z\in {\frkp }=T_oM$, the vector field $Y(t):=L_{\exp (tX)*}(Y)$ and $Z(t)=L_{\exp(tX)*}(Z)$ are respectively the parallel translation of $Y$ and $Z$ along the geodesic $\gamma(t):=_{\exp (tX)}\!\!K$ (see, for example, \cite{hel01}). Then we have
\begin{align}
\dd{}{t}((L_{\exp (tX)}^*A)(Y,Z))
=&\dd{}{t}(A_{\exp (tX)K}(L_{\exp (tX)*}(Y),
L_{\exp (tX)*}(Z)))\nnm\\
=&(\hat\nabla_{\gamma'(t)}A)(Y,Z)=0.
\end{align}
It follows that the function
\be\label{2.19-0}
(L_{\exp (tX)}^*A)(Y,Z)
=A_{\exp (tX)K}(L_{\exp (tX)*}(Y),
L_{\exp (tX)*}(Z))
\ee
is constant with respect to parameter $t$ and thus $A$ is $G$-invariant.

Conversely, we suppose that $A$ is $G$-invariant. Then for any $X,Y,Z\in {\frkp }=T_oM$, function
\eqref{2.19-0}
is again a constant along the geodesic $\gamma(t)$. Therefore,
$$
(\hat\nabla_XA)(Y,Z)=\left.\dd{}{t}\right|_{t=0} A_{\gamma(t)}(Y(t),Z(t))=0,
$$
where we have once again used the fact that $Y(t)$ and $Z(t)$ are parallel along the geodesic $\gamma(t)$. Thus $(\hat\nabla A)_o=0$. This together with the invariance of $A$ and the metric $g$ easily proves that $\hat\nabla A=0$ everywhere.
\endproof

An application of Proposition \ref{prop3.6} and \eqref{abl1} gives directly the following corollary:

{\cor\label{cor3.7} A nondegenerate equiaffine symmetric hypersurface $x:M^n\to\bbr^{n+1}$ is necessarily an affine hypersphere. In particular, the affine Gauss equation \eqref{gaus_af sph} holds in this case.}

Now we define a multiplication $\circ$ on $V$ as follows (cf. \eqref{eqn3.1} and \eqref{eqn3.2}):
\be\label{JorProd}
u\circ v:=T_u(v)\equiv T_v(u),\quad \forall u,v\in V.\ee

\begin{lem}\label{lem3.9}
$(V,\circ)$ is a nondegenerate Jordan algebra with $e$ as a unity.
\end{lem}

\proof
The multiplication $\circ$ is clearly bilinear, symmetric and has a unity $e$. Thus it needs only to verify the Jordan identity.

Arbitrarily given $u,v\in V$, we find

(1) If $u=\lambda e$ for some $\lambda\in\bbr$, then
\begin{align*}
u\circ(u^2\circ v)-u^2\circ(u\circ v)=&\lambda e((\lambda e)^2\circ v) -(\lambda e)^2\circ(\lambda e\circ v)\\
=&\lambda e\circ(\lambda^2 e\circ v)-\lambda^2 e\circ(\lambda v)=\lambda^3v-\lambda^3v=0;
\end{align*}

(2) If $v=\mu e$ for some $\mu\in\bbr$, then
\begin{align*}
u\circ(u^2\circ v)-u^2\circ(u\circ v)=&u\circ(u^2\circ \mu e) -u^2\circ(u\circ \mu e)\\
=&\mu u\circ(u^2)-\mu u^2\circ u =0;
\end{align*}

(3) If $u=X,v=Y\in V_0$, then
\begin{align*}
u\circ(u^2\circ v)=&X\circ(X^2\circ Y)=X\circ ((A(X,X)-L_1g(X,X)e)\circ Y)\\ =&X\circ (A(X,X)\circ Y-L_1g(X,X)Y)\\
=&X\circ(A(A(X,X),Y)-L_1g(A(X,X),Y)e)-L_1g(X,X)X\circ Y\\
=&A(X,A(A(X,X),Y))-L_1g(X,A(A(X,X),Y))e-L_1g(A(X,X),Y)X \\ &\ -L_1g(X,X)(A(X,Y)-L_1g(X,Y)e)\\
=&A_X(A_{A(X,X)}Y)-L_1g(A(X,X),A(X,Y))e-L_1g(A(X,X),Y)X\\ &\ -L_1g(X,X)A(X,Y)+L^2_1g(X,X)g(X,Y)e.
\end{align*}
Similarly
\begin{align*}
u^2\circ(u\circ v)=&X^2\circ(X\circ Y)=(A(X,X)-L_1g(X,X)e)\circ (A(X,Y)-L_1g(X,Y)e)\\
=&A(X,X)\circ A(X,Y)-L_1g(X,Y)A(X,X) -L_1g(X,X)A(X,Y)+L^2_1g(X,X)g(X,Y)e\\
=&A(A(X,X),A(X,Y))-L_1g(A(X,X),A(X,Y))e-L_1g(X,Y)A(X,X)\\ &\ -L_1g(X,X)A(X,Y)+L^2_1g(X,X)g(X,Y)e\\
=&A_{A(X,X)}(A_XY)-L_1g(A(X,X),A(X,Y))e-L_1g(X,Y)A(X,X) \\ &\ -L_1g(X,X)A(X,Y)+L^2_1g(X,X)g(X,Y)e.
\end{align*}
Thus we have by \eqref{gaus_af sph}
\begin{align*}
&u\circ(u^2\circ v)-u^2\circ(u\circ v)\\
=&A_X(A_{A(X,X)}Y)-A_{A(X,X)}(A_XY)+L_1(g(X,Y)A(X,X)-g(A(X,X),Y)X)\\
=&L_1(g(X,Y)A(X,X)-g(A(X,X),Y)X)-[A_{A(X,X)},A_X](Y)\\
=&R(A(X,X),X)Y\in V_0.
\end{align*}

On the other hand, for any $Z\in V_0$,
\begin{align*}
&g(u\circ(u^2\circ v)-u^2\circ(u\circ v),Z)
=g(R(A(X,X),X)Y,Z)\\
=&\fr13(g(R(Y,Z)A(X,X),X)-2g(R(Y,Z)X, A(X,X)))\\
=&\fr13(g(R(Y,Z)A(X,X),X)-2g(A(R(Y,Z)X,X),X)))\\
=&\fr13(g(R(Y,Z)A(X,X)-2A(R(Y,Z)X,X),X)\\
=&\fr13g((R(Y,Z)A)(X,X),X)=0,
\end{align*}
where the last equality is by the $G$-invariance of $A$ and Proposition \ref{prop3.6}. It then follows that
$u\circ(u^2\circ v)-u^2\circ(u\circ v)=0$ in this case.

To prove the nondegeneracy of $(V,\circ)$, let $v\in V$ and suppose that $T_v(u)=0$ for all $u\in V$. In particular, $T_v(e)=0$, that is, $v=0$.
\endproof

\begin{lem}\label{lem3.2}
Let $\lagl\cdot,\cdot\ragl$ be the inner product defined by \eqref{eqn2.18}. For all $X,Y\in V_0$, we have
\be\label{eqno3.5}
\lagl X,Y\ragl=-(n+1)L_1g(X,Y),\quad \lagl X,e\ragl=0,\quad \lagl e,e\ragl=n+1.
\ee
In particular, the inner product $\lagl\cdot,\cdot\ragl$ is nondegenerate if and only if the affine mean curvature $L_1\neq 0$, and thus in this case, $(V,\circ)$ is semi-simple as a Jordan algebra.
\end{lem}

\proof Let $\{e_i\}$ be a basis for $V_0$. Then $\{e_i,e\}$ is a basis for $V$ of which the dual basis is denoted by $\{\omega^i,\omega\}$.

Firstly for each pair $X,Y\in V_0$, we have
\begin{align*}
\lagl X,Y\ragl=&\tr T_{X\circ Y}=\omega(T_{X\circ Y}(e))+\sum_{i}\omega^i(T_{X\circ Y}(e_i))\\
=&\omega(X\circ Y)+\sum_i\omega^i((X\circ Y)\circ e_i)\\
=&\omega(A(X,Y)-L_1g(X,Y)e)+\sum_i\omega^i((A(X,Y)-L_1g(X,Y)e)\circ e_i)\\
=&-L_1g(X,Y)+\sum_i\omega^i(A(X,Y)\circ e_i-\sum_i\omega^i(L_1g(X,Y)e_i))\\
=&-L_1g(X,Y)+\sum_i\omega^i(A(A(X,Y),e_i)-L_1g(A(X,Y),e_i)e)-nL_1g(X,Y)\\
=&-(n+1)L_1g(X,Y),
\end{align*}
where $\sum_i\omega^i(A(A(X,Y),e_i))=\tr A_{A(X,Y)}=0$ due to the apolarity \eqref{apol} of $A$.

Secondly,
\begin{align*}
\lagl X,e\ragl=&\tr T_{X\circ e}=\tr T_X
=\sum_i\omega^i(X\circ e_i)+\omega(X\circ e)\\
=&\sum_i\omega^i(A(X,e_i)-L_1g(X,e_i)e)+\omega(X)\\
=&\sum_i\omega^i(A(X,e_i))=0;
\end{align*}
and
\begin{align*}
\lagl e,e\ragl=&\tr T_{e\circ e}=\tr T_e
=\sum_i\omega^i(e\circ e_i)+\omega(e\circ e)\\
=&\sum_i\omega^i(e_i)+\omega(e)
=n+1.
\end{align*}
\endproof

The above Jordan algebra $(V,\circ)$ is clearly dependent on the given point $o$ fixed at the beginning of this section. So different points on $M$ define different Jordan algebras. But due to the $G$-invariance of the Fubini-Pick form $A$, all of these Jordan algebras must isomorphic to each other:

\begin{prop}\label{prop3.1}
Let $p$ be an arbitrary point of $M$ and $(V,\td\circ)$ the Jordan algebra given by $p$. Then there exists an isomorphism $\phi:(V,\circ)\to (V,\td\circ)$ of Jordan algebras.
\end{prop}

In fact, by Proposition \ref{prop3.3}, we can write $M=G/K$ with $G$ identified with a Lie subgroup of the unimodular transformation group on $\bbr^{n+1}$ (See Remark \ref{rmk3.1}). In particular, each left action $L_g$, $g\in G$, keeps both the affine metric $g$ and the Fubini-Pick form $A$ invariant. Pick one $g\in G$ satisfying $L_g(o)=p$. Define $\phi_g:V\to V$ by
$$
\phi_g(X)=L_{g*o}(X),\quad \text{for\ }X\in V_0\equiv T_oM;\quad \phi_g(x(o))=x(p).
$$
Then it is easily checked that $\phi_g$ is an isomorphism of Jordan algebras.

\section{Jordan triples and the symmetric Lie algebras from the real Jordan algebras}\label{sec4}

By definition (\cite{nai83a}), A {\it Jordan triple} $(V,\{\cdot,\cdot,\cdot\})$ is a linear space $V$ equipped with a trilinear map $$\{\cdot,\cdot,\cdot\}:V\times V\times V\to V$$ satisfying the following two conditions:

(JT1) For all $u,v,w\in V$, $L(u,v)w=L(u,w)v$ where $L(u,v)w:=\{u,v,w\}$;

(JT2) For all $u,v,w,z\in V$, $[L(w,z),L(u,v)]=L(L(w,z)u,v)-L(u,L(w,z)v)$

\noindent
where $[\cdot,\cdot]$ is the Lie bracket for linear transformations on $V$.

A Jordan triple $(V,\{\cdot,\cdot,\cdot\})$ is called nondegenerate if the bilinear form $(\cdot,\cdot)$ on $V$ given by
\be\label{eqn4.5}
(u,v)=\tr L(u,v),\quad u,v\in V
\ee
is nondegenerate. Note that it can be verified that the bilinear form $(\cdot,\cdot)$ defined above is symmetric (\cite{nai83a}, p101).

\begin{lem}\label{lem4.1} {\rm(\cite{nai83a}), \cite{nai83b}} For a given Jordan triple $V$, we have
\be\label{eqn4.6}
(L(u,v)w,z)=(w,L(v,u)z),\quad \forall u,v,w,z\in V.
\ee
\end{lem}

Now let ${\mathcal J}=(V,\circ)$ be a nondegenerate real Jordan algebra with the unity $e$, $\dim V=n+1$. Define
\be\label{eqn4.7}
V_0=\{u\in V;\ \tr T_u=0\}.
\ee
Then it is easily seen that $V=V_0\oplus \bbr e$. In particular, when ${\mathcal J}$ is semi-simple, $V_0=e^\bot$ is the orthogonal complement of $e$ in $V$ with respect to the inner product $\lagl\cdot,\cdot\ragl$.

From the Jordan algebra $\mathcal J$, we can define a trilinear map by
\be\label{eqn4.9}
\{u,v,w\}=u\circ(v\circ w)-v\circ(u\circ w)+(u\circ v)\circ w,\quad u,v,w\in V.
\ee
Then we have
\begin{lem}\label{lem4.2} {\rm(\cite{nai83a})} $(V,\{\cdot,\cdot,\cdot\})$ given in \eqref{eqn4.9} is a Jordan triple. Furthermore, if the Jordan algebra $(V,\circ)$ is semi-simple, then $(V,\{\cdot,\cdot,\cdot\})$ is nondegenerate.
\end{lem}

We remark that, by \eqref{eqn4.9}, the two symmetric bilinear forms $(\cdot,\cdot)$ and $\lagl\cdot,\cdot\ragl$ are exactly the same.

Introduce the following linear subspaces of $\frkgl(V)$
\be\label{eqn4.10}
{\mathcal L}=\spn_\bbr\{L(u,v);\ u,v\in V\}.
\ee
Then from (JT2) it is easily seen that
\begin{lem}\label{lem4.3}(\cite{nai83a})
${\mathcal L}$ is a Lie subalgebra of $\frkgl(V)$, called the structure Lie algebra of the Jordan algebra ${\mathcal J}$.
\end{lem}

Moreover, put
\be\label{eqn4.12}
\cls=\spn_\bbr\{T_{u\circ v};\ u,v\in V\},\quad
\frkk={\rm gen}\{[T_u,T_v];\ u,v\in V\}.
\ee
Then

\begin{lem}\label{lem4.4} It holds that
\be\label{eqn4.13}
\cls=\spn_\bbr\{T_u;\ u\in V\},\quad \frkk={\rm gen}\{[T_X,T_Y];\ X,Y\in V_0\}\subset\frkgl(V).
\ee
Furthermore, $\frkk$ is a Lie subalgebra of ${\mathcal L}$ and ${\mathcal L}=\frkk\oplus \cls$.
\end{lem}

\proof (1) The first equality in \eqref{eqn4.13} is apparent since the existence of the unity $e$. To prove the second one it suffices to show that \be\label{eqn4.14}\frkk\subset{\rm gen}\{[T_X,T_Y];\ X,Y\in V_0\}.\ee
For any $X,Y\in V_0$, $\lambda,\mu\in\bbr$, and any $v\in V$,
\begin{align*}
T_{X+\lambda e}(T_{Y+\mu e}(v))=&T_{X+\lambda e}((Y+\mu e)\circ v)
=(X+\lambda e)\circ(Y\circ v+\mu v)\\
=&X\circ(Y\circ v)+\mu(X\circ v)+\lambda(Y\circ v) +\lambda\mu e\\
=&T_X(T_Y(v))+(\mu X+\lambda Y)\circ v+\lambda\mu e;\\
T_{Y+\mu e}(T_{X+\lambda e}(v))=&T_Y(T_X(v))+(\lambda Y+\mu X)\circ v+\lambda\mu e.
\end{align*}
Thus
$$
[T_{X+\lambda e},T_{Y+\mu e}](v)=T_{X+\lambda e}(T_{Y+\mu e}(v))-T_{Y+\mu e}(T_{X+\lambda e}(v)) =[T_X,T_Y](v),
$$
that is, $[T_{X+\lambda e},T_{Y+\mu e}]=[T_X,T_Y]$ which proves \eqref{eqn4.14}.

(2) For any pair of $u,v,w\in V$, we have by \eqref{eqn4.9}
\begin{align*}
L(u,v)w=&\{u,v,w\}=u\circ(v\circ w)+(u\circ v)\circ w-v\circ (u\circ w)\\
=&T_u(T_v(w))-T_v(T_u(w))+T_{u\circ v}(w)
=([T_u,T_v]+T_{u\circ v})w.
\end{align*}
Thus we obtain
\be
L(u,v)=[T_u,T_v]+T_{u\circ v}.
\ee
Consequently
\be\label{eqn4.16}
L(u,v)+L(v,u)=2T_{u\circ v},\quad L(u,v)-L(v,u)=2[T_u,T_v].
\ee
Clearly, the decomposition ${\mathcal L}=\frkk\oplus \cls$ comes right from \eqref{eqn4.16} and (1).

(3) The second equality of \eqref{eqn4.16} indicates that $\frkk\subset\cll$ as a linear subspace, which together with the fact that $\cll$ is also a Lie subalgebra of $\frkgl(V)$ shows that $\frkk$ is a Lie subalgebra of $\cll$. This can also be shown in another more direct manner: First note that by (JT2)
\be\label{eqn4.18}
[\Phi,L(u,v)]=L(\Phi(u),v)-L(u,\Phi^t(v)),\quad \forall \Phi\in {\mathcal L},\,\forall u,v\in V.
\ee

Since $\Phi^t=-\Phi$ for any $\Phi\in \frkk$, it follows that for $u,v\in V$
\begin{align*}
[\Phi,[T_u,T_v]]=&\fr12[\Phi,L(u,v)-L(v,u)]=\fr12([\Phi,L(u,v)]-[\Phi,L(v,u)])\\
=&\fr12(L(\Phi(u),v)-L(u,\Phi^t(v))-L(\Phi(v),u)+L(v,\Phi^t(u))) \\
=&\fr12(L(\Phi(u),v)+L(u,\Phi(v))-L(\Phi(v),u)-L(v,\Phi(u)))\\
=&\fr12(L(\Phi(u),v)-L(v,\Phi(u))+L(u,\Phi(v))-L(\Phi(v),u)))\\
=&[T_{\Phi(u)},T_v]+[T_u,T_{\Phi(v)}]\in \frkk,
\end{align*}
implying that the bracket $[\cdot,\cdot]$ of $\cll$ is closed on $\frkk$.
\endproof

We claim that, if the Jordan algebra ${\mathcal J}$ is semi-simple, then
\be\label{eqn4.19}
\Phi(V)\subset V_0,\text{\ for all\ }\Phi\in \frkk.
\ee
In fact, for any $\Phi\in \frkk$, $\Phi^t=-\Phi$. By the second equality of \eqref{eqn4.13} we can assume without loss of generality that $\Phi=[T_X,T_Y]$ for some $X,Y\in V_0$. It follows that for any $u\in V$,
\be\label{eqn4.20}
\lagl\Phi(v),e\ragl=\lagl[T_X,T_Y](v),e\ragl=-\lagl v,[T_X,T_Y](e)\ragl=0
\ee
implying that $\Phi(u)\in e^\bot\equiv V_0$.

{\prop For any $\Phi\in \frkk$, we have
\be\label{eqn4.21}
[\Phi,T_u]=T_{\Phi(u)},\quad\forall u\in V.
\ee}

\proof For any $\Phi\in \frkk$, $u\in V$ we find by \eqref{eqn4.18}
\begin{align*}
[\Phi,T_u]=&\fr12[\Phi,L(u,e)+L(e,u)]\\
=&\fr12(L(\Phi(u),e)-L(u,\Phi^t(e))) +\fr12(L(\Phi(e),u)-L(e,\Phi^t(u)))\\
=&\fr12(L(\Phi(u),e)+L(u,\Phi(e))) +\fr12(L(\Phi(e),u)+L(e,\Phi(u)))\\
=&\fr12(L(\Phi(u),e)+L(e,\Phi(u)))\\
=&T_{\Phi(u)\circ e}
=T_{\Phi(u)}.
\end{align*}
\endproof

From \eqref{eqn4.21} the following corollary is easily checked:

{\cor\label{cor4.6} The Jordan product $\circ$ is invariant by $\frkk$, that is
\be\label{eqn4.22}
\Phi(u\circ v)=\Phi(u)\circ v+u\circ\Phi(v),\quad \forall \Phi\in\frkk,\ \forall u,v\in V.
\ee}

From the Lie algebra ${\mathcal L}$ we can define a new Lie algebra. In fact, if we put
\be\label{eqn4.23}
\frkp=\{T_X;\ X\in V_0\},\quad \frkg=\frkk\oplus \frkp,
\ee
then $\frkg$ is a Lie subalgebra of the structure Lie algebra ${\mathcal L}$ since $[e,u]=0$ for all $u\in V_0$.

{\dfn\rm The Lie algebra $\frkg$ given in \eqref{eqn4.23} is called the the restricted structure Lie algebra of the Jordan algebra ${\mathcal J}=(V,\circ)$.

Define $\theta:\frkg\to \frkg$ by $\theta(\Phi+T_X)=\Phi-T_X$ for all $\Phi\in\frkk$, $X\in V_0$. Then we have

\begin{prop}\label{prop4.1} $({\mathfrak g},\theta)$ is an effective symmetric Lie algebra, or equivalently, $({\mathfrak g},\frkk)$ is an effective symmetric pair of Lie algebras.
\end{prop}

\proof We have show that both $\frkg$ and $\frkk$ are Lie subalgebras of $\frkgl(V)$. Then the second equality of \eqref{eqn4.13} and \eqref{eqn4.14} indicate that $(\frkg,\frkk)$ is a symmetric pair of Lie algebras; The effectiveness of $(\frkg,\frkk)$ follows from the proof of Lemma \ref{lem3.3}.\endproof

\section{Nondegenerate hypersurfaces associated with semi-simple real Jordan algebras}\label{sec5}

Let ${\mathcal J}=(V,\circ)$ be a semi-simple real Jordan algebra and $(\frkg,\theta)$ the symmetric Lie algebra given in the Last section where $\frkg$ is the restricted Lie algebra of ${\mathcal J}$. By Proposition \ref{prop2.3}, ${\mathcal J}$ is unital. Let $e$ be a unity of ${\mathcal J}$. In this section we shall use ${\mathcal J}$ to construct a connected, simply connected and nondegenerate equiaffine hypersurface with a given nonzero affine mean curvature $L_1$. To this end we first need to find the symmetric pair $(G,K)$ of Lie groups associated with the symmetric Lie algebra $(\frkg,\theta)$, where $G$ should be connected and simply connected and $K$ should be connected. We begin with introducing the group
\be\label{eqn5.1}
\bar G={\rm Gen}\{P_u;\ u\in V,\ \det(P_u)=1\}\subset SL(V),
\ee
where we use ``${\rm Gen}\{\cdots\}$'' to denote the Lie group generated by a subset ``$\{\cdots\}$'' and $SL(V)$ is the special linear group on $V$. By virtue of Proposition \ref{prop2.6}, it is easy to verify that $\bar G$ is a Lie subgroup of $SL(V)$ which preserves the determinant defined on $V$ (see \eqref{eqn2.17}). Define $G$ to be the universal covering group of the unit component $\bar G_0$ of $\bar G$ and denote by $\pi:G\to \bar G_0$ the covering homomorphism.

{\prop\label{prop5.1} The Lie algebra of $G$ is precisely the restricted structure Lie subalgebra $\frkg$ introduced previously in Section \ref{sec4}.}

\proof Denote by $\cll_G$ be the Lie algebra of $G$ on $V$. Then $\cll_G$ equals to the Lie algebra $\cll_{\bar G_0}$ of the unit component $\bar G_0$ of $\bar G$ and therefore equals to the Lie algebra $\cll_{\bar G}$ of $\bar G$. Thus it follows that
$$
\cll_G=\cll_{\bar G}={\rm gen}\{\td\gamma'(0);\ \td\gamma:(-\veps,\veps)\to \bar G\in C^\infty,\ \td\gamma(0)=I_V,\ \veps>0\}\subset \frksl(V),
$$
where $I_V$ is the identity map in the general linear group $GL(V)$. For any smooth curve $\td\gamma(t)$ ($-\veps<t<\veps$) around the identity in $G$, we can assume that there is a curve $u=u(t)$, nondegenerate for each $t\in (-\veps,\veps)$, with $u(0)=e$ such that $\td\gamma(t)=P_{u(t)}$. Then it follows that
\begin{align}
\td\gamma'(0)=&\left.\dd{}{t}\right|_{t=0}(P_{u(t)}) =\left.\dd{}{t}\right|_{t=0}(2T^2_{u(t)}-T_{u^2(t)})\nnm\\
=&4T_{u'(0)}-2T_{u'(0)}=2T_{u'(0)}\in \cls\label{eqn5.2}.
\end{align}

On the other hand, since $\det P_{u(t)}\equiv 1$ it follows that
\be\label{eqn5.3}
0=\left.\dd{}{t}\right|_{t=0}(\det P_{u(t)})=\tr\left( \left.\dd{}{t}\right|_{t=0} P_{u(t)}\right)=2\tr(T_{u'(0)}).
\ee
Consequently we have $\td\gamma'(0)=T_{u'(0)}\in \frkp$. Thus
\be\label{eqn5.4}
\cll_G={\rm gen}\{\td\gamma'(0)\}\subset {\rm gen}(\frkp)=\frkg.
\ee

Conversely, for any $X\in V_0$, let
$$u(t)=e+\fr12tX \text{\ and\ } \td u(t)=\left(\det P_{u(t)}\right)^{-\fr1{2(n+1)}}u(t),\quad -\veps<t<\veps.
$$
Sine $u(0)=e$, $u'(0)=\fr12X$, it follows that
$P_{u(0)}=I_V$, $\td u(0)=u(0)=e$
and, by \eqref{eqn5.2} and \eqref{eqn5.3}, that
$$
\left.\dd{}{t}\right|_{t=0}P_{u(t)}=2u'(0)=X,\quad \left.\dd{}{t}\right|_{t=0}\det P_{u(t)}=2\tr T_{u'(0)}=\tr T_X=0.
$$
Therefore
\begin{align}
\td u'(0)=&\left.\dd{}{t}\right|_{t=0}\left((\det P_{u(t)})^{-\fr1{2(n+1)}}u(t)\right)\nnm\\
=&-\fr1{2(n+1)}(\det P_{u(0)})^{-\fr{2n+3}{2(n+1)}} \left.\dd{}{t}\right|_{t=0}(\det P_{u(t)})u(0)+\left((\det P_{u(0)})^{-\fr1{2(n+1)}}\right)u'(0)\nnm\\
=&\fr12X.\label{eqn5.5}
\end{align}
Define
$\td\gamma(t)=P_{\td u(t)}$, $-\veps<t<\veps$.
Then
\begin{align*}
\det\td\gamma(t)=&\det P_{\td u(t)}=\det\left(P_{(\det P_{\td u(t)})^{-\fr1{2(n+1)}}u(t)}\right)\\
=&\det\left((\det P_{\td u(t)})^{-\fr1{n+1}}P_{u(t)}\right)
=(\det P_{\td u(t)})^{-1}\det P_{\td u(t)}=1
\end{align*}
implying that $\td\gamma(t)\in G$ for any $t$. It follows that
\be
T_X=2T_{\td u'(0)}=\left.\dd{}{t}\right|_{t=0}\left(P_{\td u(t)}\right)=\td\gamma'(0)
\in \cll_G,
\ee
showing that $\frkg\subset\cll_G$. This with \eqref{eqn5.4} proves that $\cll_G=\frkg$. \endproof

Next we define $K={\rm Gen}\{\exp_G(X);\ X\in \frkk\}$. Then $K$ is a connected and closed Lie subgroup of $G$ consisting of all element $g$ in $G$ satisfying $\pi(g)(e)=e$, where $\pi:G\to \bar G_0$ is the covering homomorphism. Since $G$ is simply connected, the involutive automorphism $\theta:\frkg\to\frkg$ lifts to an automorphism $\sigma:G\to G$. Therefore $(G,K)$ is a symmetric pair and thus $M^n:=G/K$ is a connected and simply connected symmetric space. Denote by $o={}_{I_V}K\in M^n$. Then the tangent space $T_oM^n$ is identified with $\frkp=\{T_X;\ X\in V_0\}$ which is in turn identified with $V_0$.

For a given constant $L_1\neq 0$, choose  $$C=-\sgn(L_1)\sqrt{n+1}((n+1)|L_1|)^{-\fr12(n+2)},$$
and define a map $f:G\to V\equiv\bbr^{n+1}$ by $f(g)=C\pi(g)(e)$, $g\in G$. Clearly, for any $g_1,g_2\in G$,
$f(g_1)=f(g_2)$ if and only if $\pi(g_1^{-1}\circ g_2)(e)=e$, that is, $g_1^{-1}\circ g_2\in K$. Therefore, $f$ naturally induces a smooth map $x:M^n\to V\equiv\bbr^{n+1}$ by
\be\label{eqn5.7}
x(gK)=C\pi(g)(e),\quad\forall g\in G.
\ee

{\prop\label{prop5.2} $x$ is a connected, simply connected and nondegenerate equiaffine symmetric hypersurface with the given nonzero constant $L_1$ as its affine mean curvature.}

\proof
Since $\pi(G)\subset SL(V)$, we can choose a volume element on $\bbr^{n+1}$,
say, the canonical volume element with respect to the inner product $\lagl\cdot,\cdot\ragl$ on $V$, so that $\pi(G)$ can be identified with a subgroup of the
group ${\rm UA}(n+1)$ of unimodular affine transformation on $\bbr^{n+1}$.
Therefore, the induced map $x$ is equivariant as an affine hypersurface
in $\bbr^{n+1}$. Consequently all the equiaffine invariants of $x$ such as the affine metric, the Fubini-Pick form and the affine second fundamental form are ${\rm G}$-invariant.

Now for each $T_X\in\frkp\equiv T_oM^n$, $X\in V_0$, $a(t):={}_{\exp tT_X}K$ is a geodesic curve on $M^n$ with respect to any $G$-invariant psedu-Riemannian metric. It holds clearly that
$$
x_*(T_X)=\left.\dd{}{t}\right|_{t=0}x(a(t))=C\left.\dd{}{t}\right|_{t=0}(\exp tT_ X(e))=CT_X(e)=C(X\circ e)=C\cdot X,
$$
where the isomorphism $\pi_*$ is omitted. This shows that $x$ is an immersion at $o$ and thus is an immersion globally since $x$ is equivariant. Clearly, $x$ is injective and is thus an imbedding of $M^n$ into $\bbr^{n+1}$.

Moreover, since for each $X\in V_0$,
$$\lagl X,e\ragl=\tr(X\circ e)=\tr X=0,\quad \lagl e,e\ragl=n+1,$$
$x(o)$ is a transversal vector of $x$ at $o$ and thus is transversal everywhere. Furthermore,
for an arbitrary $Y\in V_0$, denote by $Y^*$ the Killing vector field on $M^n$ induced by $T_Y$, then the value of $Y^*$ at $a(t)$
$$
Y^*|_{a(t)}=\left.\dd{}{s}\right|_{s=0}(\,{}_{\exp sT_Y a(t)}K)=\left.\dd{}{s}\right|_{s=0}({}\,_{\exp sT_Y \exp tT_X}K).
$$
Therefore
$$
x_*(Y^*|_{a(t)})=C\dd{}{t}\left.\dd{}{s}\right|_{s=0}(\exp sT_Y \exp tT_X(e)).
$$
It follows that
\begin{align}
T_X(x_*(Y^*))=&C\left.\ppp{}{t}{s}\right|_{t=s=0}(\exp s T_Y\exp tT_X(e))\nnm\\
=&C(Y\circ(X\circ e))=C(Y\circ X)\nnm\\
=&C\left(X\circ Y-\fr1{n+1}\tr (T_{X\circ Y})e\right)+\fr C{n+1}\lagl X,Y\ragl e.\label{gaussf}
\end{align}
Since
$$
\lagl X\circ Y-\fr1{n+1}\tr (T_{X\circ Y})e,e\ragl =\tr T_{X\circ Y}-\fr1{n+1}\tr(T_{X\circ Y})\lagl e,e\ragl=0,
$$
we have
$$
C(X\circ Y-\fr1{n+1}\tr (T_{X\circ Y})e)\in x_*(TM^n).
$$
This with \eqref{gaussf} implies that $x$ is nondegenerate since $\lagl X,Y\ragl =\tr T_{X\circ Y}$ is.

Note that, by Corollary \ref{cor4.6}, the inner product $\lagl\cdot,\cdot\ragl$ on $V_0$ is $K$-invariant. It follows that the affine metric $g$ of $x$ is exactly the invariant metric on $G/K$ induced by
$$
g_o(T_X,T_Y)\equiv g_o(X,Y):=\sgn(C)\left(\fr {|C|}{\sqrt{n+1}}\right)^{\fr2{n+2}}\lagl X,Y\ragl=-\fr1{(n+1)L_1}\lagl X,Y\ragl,\quad \forall\, X,Y\in V_0
$$
(see \eqref{affine metric}). Clearly, $g$ is a symmetric metric since $g_o$ is invariant by the involution $\theta(T_ X)=-T_X$, $X\in V_0$.

Let $A_o$ be the $(1,2)$ tensor on $V_0\equiv\frkp$ defined by
\be\label{eqn5.9}
A_o(X,Y)=X\circ Y-\fr1{n+1}\tr(T_{X\circ Y})e,\quad \forall X,Y\in V_0,
\ee
which gives a linear map for any $X\in  V_0$: $A_o(X): V_0\to  V_0$ by $A_o(X)Y=A_o(X,Y)$, $Y\in  V_0$.

Define
\be\label{eqn5.10}
A_o(X,Y,Z)=g_o(A_o(X,Y),Z),\quad \forall X,Y,Z\in V_0.
\ee
Then it is necessary for us to make the following two claims:

{\bf Claim (1)} The $(0,3)$-tensor $A_o(X,Y,Z)$ is totally symmetric for all $X,Y,Z\in V_0$.

In fact, the symmetry of $X,Y$ is obvious; Moreover, by the definition of $A_o$,
\begin{align*}
A_o(X,Y,Z)=&g_o(A_o(X,Y),Z)=-\fr1{(n+1)L_1}\lagl A_o(X,Y),Z\ragl\\
=&-\fr1{(n+1)L_1}(\lagl X\circ Y-\fr1{n+1}\tr(T_{X\circ Y})e,Z\ragl\\
=&-\fr1{(n+1)L_1}(\lagl T_X(Y),Z\ragl -\fr1{(n+1)^2L_1}\tr(T_{X\circ Y})\lagl e,Z\ragl\\ =&-\fr1{(n+1)L_1}(\lagl T_X(Z),Y\ragl -\fr1{(n+1)^2L_1}\tr(T_{X\circ Z})\lagl e,Y\ragl\\
=&A_o(X,Z,Y),\end{align*}
where we have used the symmetry of $T_X$ and the fact that $V_0$ is the orthogonal complement of $e$. Claim (1) is then proved.

{\bf Claim (2)} For each $X\in  V_0$, the linear map $A_o(X)$ is traceless.

In fact, if $\{e_i\}$ is a basis for the linear space $V_0$, then $\{e_i,e\}$ is a basis for $V$ of which the dual basis is denoted by $\{\omega^i,\omega\}$. It follows that
\begin{align*}
\tr A_o(X)=&\sum\omega^i(A_o(X)e_i)=\sum\omega^i(A_o(X,e_i))\\
=&\sum\omega^i(X\circ e_i-\fr1{n+1}\tr(T_{X\circ e_i})e)\\
=&\sum\omega^i(X\circ e_i)+\omega(T_X(e))\\
=&\tr(T_X)=0,
\end{align*}
proving Claim (2).

Now, since $Y^*$ is chosen to be the Killing vector field on $M$ corresponding to $Y$, we have $\hat\nabla_XY^*=0$ where $\hat\nabla$ is the Levi-Civita connection of the affine metric $g$. It follows that the symmetric three form $A_o$ defined by \eqref{eqn5.10} coincides with the Fubini-Pick form of $x$ at the origin $o$.

Taking the trace of \eqref{gaussf} with respect to $g_o$ and using Claims (1) and (2) given above, we find that, at $o$, the affine normal vector
$$
\xi_o=\fr1n(\Delta_g x)_o=-(n+1)L_1\cdot\fr{C}{n+1}\cdot e=-L_1x(o).
$$
$\xi_o$ is clearly invariant by $K$ and for any $X\in \frkg$, $X(\xi_o)\in V_0=x_*(T_o(M))$. Then the equivariant transversal vector field $\xi$ induced by $\xi_o$ coincides with the affine normal vector (see Lemma 3.4 in \cite{nom-sas94b}). Since $x$ is also equivariant, $\xi=-L_1x$ holds identically, which implies that $x$ is a proper affine hypersphere with the origin as its affine center and $L_1$ as its affine mean curvature.

Finally, the equivariance of $x$ and the affine uniqueness theorem (Theorem \ref{af uniq}) imply that its Fubini-Pick form $A$ is $G$-invariant, which is uniquely determined by the cubic form $A_o$ given by \eqref{eqn5.9} and \eqref{eqn5.10}. In particular, $x$ is an affine symmetric hypersphere in $\bbr^{n+1}$.
\endproof

{\prop\label{prop5.3} The image $x(M^n)$ in $\bbr^{n+1}$ of the embedding $x$ defined by \eqref{eqn5.7} is identical to the connected component $\bar M^n$ containing the vector $Ce$ of the following hypersurface
$$
\td M^n=\{u\in V;\ \det (P_u)=C^{2(n+1)}\}.
$$}

\proof Firstly note that $G$ is determinant-preserving, $\det P_{Ce}=\det(C^2I_V)=C^{2(n+1)}$ and $M^n$ is connected. These things indicate that the image $x(M^n)$ is contained in $\bar M^n$. On the other hand, the fact that $\dim x(M^n)=\dim \bar M^n$ and the completeness make true the inverse inclusion $\bar M^n\subset x(M^n)$.
\endproof

\section{The proof of Theorem \ref{main}}\label{sec6}

In this section, we shall sum up discussions in Sections \ref{sec3}, \ref{sec4} and \ref{sec5} to complete the proof of our first main theorem, that is, Theorem \ref{main}

To begin the proof, let $L_1\neq 0$ be arbitrarily given and let $x:M^n\to \bbr^{n+1}$ be a connected, simply connected and nondegenerate equiaffine symmetric hypersurface with the affine metric $g$ and affine mean curvature $L_1$. Then $(M^n,g)$ is by definition a connected pseudo-Riemannian symmetric space $G/K$ and the Fubini-Pick form $A$ is $G$-invariant, where the groups $G$ and $K$ are given in Section \ref{sec3} (see \eqref{eqn3.3}, \eqref{eqn3.4} and Proposition \ref{prop3.3}).

Choose a constant $C$ given by \eqref{C} and put $o=eK\in G/K$, $V_0=T_o(G/K)$, $V=V_0\oplus\bbr\cdot e$ with $e=C^{-1}x(o)$. Then by Lemma \ref{lem3.9}, the Jordan algebra ${\mathcal J}\equiv (V,\circ)$ constructed from $x$ is a semi-simple real Jordan algebra of dimension $n+1$ with $e$ its unity where the Jordan product $\circ$ is defined via \eqref{eqn3.1}, \eqref{eqn3.2} and \eqref{JorProd}. Note that by Proposition \ref{prop3.1}, different choices of the origin $o$ give different but isomorphic Jordan algebras.

Conversely, given a semi-simple real Jordan algebra ${\mathcal J}=(V,\circ)$ of dimension $n+1$ with a unity $e$, let $V_0=e^\bot$ be the orthogonal complement of the unity $e$ with respect to the nondegenerate product $\lagl\cdot,\cdot\ragl$ defined by \eqref{eqn2.18}, which is an $n$-dimensional subspace of $V$. Then we have a nondegenerate Jordan triple $(V,\{\cdot,\cdot,\cdot\})$, (see Lemma \ref{lem4.2}) defined by \eqref{eqn4.9}, from which the structure Lie algebra $\cll$ is defined (see \eqref{eqn4.10}, \eqref{eqn4.12}, Lemma \ref{lem4.3} and Lemma \ref{lem4.4}) which is a subalgebra of the special linear Lie algebra $\frksl(V)$. The Lie algebra $\cll$ has a Lie subalgebra $g$, that is, the restricted structure Lie algebra (see \eqref{eqn4.23}) containing a Lie subalgebra $\frkk$ defined in \eqref{eqn4.13} which keeps the Jordan product $\circ$ invariant (see Corollary \ref{cor4.6}) and, by Lemma \ref{lem4.1}, $(\frkg,\frkk)$ is an effective pair of symmetric Lie algebras.

Let $G$ be the universal covering group of the unit component $\bar G_0$ of $\bar G$ given by \eqref{eqn5.1}, and $K$ be defined by $K={\rm Gen}\{\exp_G(X);\ X\in \frkk\}$. Then $K$ is a connected and closed Lie subgroup of $G$ which leaves the unity $e$ invariant, and $M^n:=G/K$ is a connected and simply connected symmetric space. Then by Proposition \ref{prop5.2}, the map $x:M^n\to V\equiv \bbr^{n+1}$ defined by \eqref{eqn5.7} is a connected, simply connected and nondegenerate equiaffine symmetric hypersurface with the given constant $L_1$ as its affine mean curvature.

Furthermore, to make clear the one-to-one correspondence stated in the theorem, we should clarify that, if the semi-simple Jordan algebra ${\mathcal J}\equiv (V,\circ)$ is derived, as described in Section \ref{sec3}, from a given connected, simply connected and nondegenerate hypersurface $x:M^n\to\bbr^{n+1}$ with the affine mean curvature $L_1\neq 0$, and $\td x:\td G/\td K\to V\equiv\bbr^{n+1}$ is the connected, simply connected and nondegenerate hypersurface defined in \eqref{eqn5.7} where $(\td G,\td K)$ is the symmetric pair associated with the symmetric pair $(\td\frkg,\td\frkk)$ with $\td\frkg=\td\frkk\oplus\td\frkp$ being the restricted structure Lie algebra of ${\mathcal J}$, then $M^n\cong \td G/\td K$ and $\td x$ is affine equivalent to $x$. For doing this, let $M^n\equiv G/K$ with $(G,K)$ the symmetric pair of Lie groups determined in Proposition \ref{prop3.3}, the corresponding pair $(\frkg,\frkk)$ being defined in \eqref{eqn3.3}, \eqref{eqn3.3-1} and \eqref{eqn3.4} (cf. Lemma \ref{lem3.3}).

{\lem\label{lem6.1} These two symmetric Lie algebra pairs $(\frkg,\frkk)$ and $(\td\frkg,\td\frkk)$ are precisely the same.}

{\it Proof of Lemma \ref{lem6.1}}:

In fact, we only need to show that Lie algebras $\frkk$ and $\td\frkk$ coincide to each other and the Lie bracket given by \eqref{eqn3.4} and that of $\cll\subset \frkgl(V)$ are the same. To this end, we use Corollary \ref{cor3.7} and the affine Gauss equation \eqref{gaus_af sph} to compute, for all $X,Y,Z\in V_0\equiv T_o(G/K)$,
\begin{align*}
[T_X,T_Y](Z)=&(T_X\circ T_Y-T_Y\circ T_X)(Z)=T_X(T_Y(Z))-T_Y(T_x(Z))\\
=&(X\circ(Y\circ Z))-(Y\circ(X\circ Z))\\
=&X\circ (A(Y,Z)-L_1g(Y,Z)e)-Y\circ (A(X,Z)-L_1g(X,Z)e)\\
=&A(X,A(Y,Z))-L_1g(X,A(Y,Z))e-L_1g(Y,Z)X\\
&-A(Y,A(X,Z))+L_1g(Y,A(X,Z))e+L_1g(X,Z)Y\\
=&-(L_1(g(Y,Z)X-g(X,Z)Y)-A_X(A_Y(Z))+A_Y(A_X(Z)))\\
=&-(L_1(g(Y,Z)X-g(X,Z)Y)-[A_X,A_Y](Z))=-R(X,Y)(Z),
\end{align*}
that is,
\be\label{eqn6.1}[T_X,T_Y]=-R(X,Y),\quad \forall X,Y\in V_0\ee
implying that the two Lie algebras $\frkk,\ \td\frkk$ are the same. Then from \eqref{eqn6.1} together with \eqref{eqn4.21} easily follows the lemma.

Since $G$ and $\td G$ are connected and simply connected, Lemma \ref{lem6.1} indicates that the symmetric pairs $(G,K)$ and $(\td G,\td K)$ mentioned above are isomorphic to each other. Thus we can assume that $\td G=G$, $\td K=K$ implying that $M^n=G/K=\td G/\td K$. On the other hand, by the affine uniqueness theorem, $x$ is also $G$-equivariant. It follows that
$$
\td x(gK)=Cg\cdot e=Cg\cdot C^{-1}x(o)=gx(o)=x(gK),\quad \forall g\in G.
$$

Finally, if two equiaffine symmetric hypersurfaces $x$ and $\td x$ are affine equivalent, then they have the same affine metric $g$, the same Fubini-Pick form $A$ and the same affine mean curvature $L_1$. Consequently, from the definition of the Jordan product (see \eqref{eqn3.1},\eqref{eqn3.2} and \eqref{JorProd}), it follows that $x$ and $\td x$ correspond to isomorphic semi-simple Jordan algebras. Conversely, isomorphic Jordan algebras certainly have isomorphic structure Lie algebras, and thus the corresponding Jordan algebra isomorphism restricts to isomorphic subspaces that are respectively the orthogonal complements of the two unities. It follows that the resulting restricted structure Lie algebras are isometric. This mean that the associated effective symmetric pairs of Lie algebras by \eqref{eqn4.23} are isomorphic to each other, which in turn determine isomorphic connected and simply connected symmetric pair of groups. Then by the construction of last section, they will give isomorphic symmetric spaces and the resulting equiaffine hypersurfaces given by \eqref{eqn5.7} are clearly affine equivalent. This finishes the proof of Theorem \ref{main}.

\section{The Calabi-type composition of equiaffine symmetric nondegenerate hypersurfaces}\label{sec7}

For the need of the final classification theorem, we introduce and discuss a Calabi-type composition of nondegenerate equiaffine symmetric hypersurfaces.

Let $x_\alpha:M^{n_\alpha}_\alpha\to V_\alpha$, $n_\alpha\neq 0$, $\alpha=1,\cdots,r$ ($r\geq 2$), be nondegenerate equiaffine hypersurfaces with affine mean curvature $\lalp\neq 0$. Without loss of generality we can assume that all $M^{n_\alpha}_\alpha$'s are connected and simply connected. Then by Theorem \ref{main}, each $x_\alpha$ defines a real semi-simple Jordan algebra ${\mathcal J}_\alpha=(V_\alpha,\circ)$ with unity $e_\alpha:=C^{-1}_\alpha x(o_\alpha)$ up to isomorphism where $C_\alpha$ is determined by $\lalp$ via \eqref{C}, $\alpha=1,\cdots,r$. Define $V=V_1\oplus\cdots\oplus V_r$. Then $V$ can be made into a Jordan algebra ${\mathcal J}:=(V,\circ)$ such that ${\mathcal J}_\alpha$'s are Jordan algebra ideas of ${\mathcal J}$, that is
$$
u\circ v:=u_1\circ v_1+\cdots+u_r\circ v_r,\quad u=\sum u_\alpha,\ v=\sum v_\alpha, \quad u_\alpha,v_\alpha\in V_\alpha,\ 1\leq\alpha\leq r.
$$
Let $e=e_1+\cdots+e_r$. Then it is easy to check that ${\mathcal J}$ is also unital with unity $e$.

For any $u=u_1+\cdots+u_r$, $v=v_1+\cdots+v_r$, $u_\alpha,v_\alpha\in V_\alpha$ ($1\leq\alpha\leq r$), the symmetric bilinear form $\lagl u,v\ragl_{\mathcal J}$ is expressed by
$$
\lagl u,v\ragl_{\mathcal J}=\tr\!_{\mathcal J}(T_{u\circ v})=\sum_\alpha\tr\!_{{\mathcal J}_\alpha}(T_{u_\alpha\circ v_\alpha}) =\sum_\alpha\lagl u_\alpha,v_\alpha\ragl_{{\mathcal J}_\alpha}.
$$
It follows that ${\mathcal J}$ is semi-simple and the orthogonal complement $V_0=e^\perp$ of the unity $e$ in $V$ is decomposed by
$$
V_0=(\oplus_\alpha V_{\alpha0})\oplus \bar\frkp_0,
$$
where $V_{\alpha0}\equiv e^\perp_\alpha$ is the orthogonal complement of the unity $e_\alpha$ in $V_\alpha$, and
$$
\bar\frkp_0=\{t_1e_1+\cdots+t_re_r;\ t=(t_1,\cdots,t_r)\in\bbr^r,\ \sum_\alpha(n_\alpha+1)t_\alpha=0\}.
$$

Clearly, if we put $\frkp_0=\{T_X;\ X\in\bar\frkp_0\}$, then $\frkp_0$ is a commutative Lie subalgebra of $\frkg$ and contained in the center ${\mathfrak z}(\frkg)$ of $\frkg$, which can be identified with the $(r-1)$-dimensional Lie subalgebra $T_0$ of the commutative Lie algebra $\bbr^r$:
$$
T_0:=\{t=(t_1,\cdots,t_r)\in\bbr^r;\ (n_1+1)t_1+\cdots+(n_r+1)t_r=0\}.
$$
The action of $\frkp_0$ on $V$ is component-wise and thus is identical to that of the Lie subalgebra $T_0\subset\bbr^r$ on $V$ which is also component-wise, that is, for each $t=(t_1,\cdots,t_r)\in T_0$ and $u=u_1+\cdots+u_r\in V$ with $u_\alpha\in V_\alpha$, $1\leq \alpha\leq r$,
\be
t\cdot u=t_1u_1+\cdots+t_ru_r.
\ee

For $\alpha=1,\cdots,r$, let $\frkg_\alpha=\frkk_\alpha\oplus \frkp_\alpha$ be the restricted structure Lie algebra of ${\mathcal J}_\alpha$, and $(G_\alpha,K_\alpha)$ the symmetric pair of Lie groups corresponding to the symmetric pair $(\frkg_\alpha,\frkk_\alpha)$ of Lie algebras, with $G_\alpha$ being connected, simply connected and $K$ connected. Then $\frkk_\alpha=[\frkp_\alpha,\frkp_\alpha]$, and $M^{n_\alpha}_\alpha=G_\alpha/K_\alpha$ as symmetric spaces. Moreover, the canonical decomposition of the restricted structure Lie algebra of ${\mathcal J}$ is $\frkg=\frkk\oplus\frkp$, where
\be\label{eqn7.2}
\frkp=(\oplus_\alpha\frkp_\alpha)\oplus\frkp_0,\quad\frkk= [\frkp,\frkp]=\oplus_\alpha[\frkp_\alpha,\frkp_\alpha]=\oplus_\alpha\frkk_\alpha,
\ee
since $[\frkp_\alpha,e_\alpha]=[\frkp_\alpha,\frkp_\beta]=0$ for all $\alpha$ and $\beta\neq\alpha$.

On the other hand, denote by $G\subset GL(V)$ the connected Lie group with its Lie algebra being the restricted Lie algebra $\frkg$ of ${\mathcal J}$, and $K\subset G$ be the connected Lie subgroup of $G$ with Lie subalgebra $\frkk$. Then by the proof of Proposition \ref{prop5.2}, the map
\be\label{eqn7.3}
x:G/K\to V, \text{\ given by\ } gK\mapsto x(gK):=Cg\cdot e,\quad\forall g\in G,
\ee
is a nondegenerate equiaffine symmetric hypersurface, where $C$ is the constant by \eqref{C} with an arbitrarily given $L_1\neq 0$.

Moreover, each $g_\alpha\in G_\alpha \subset GL(V_\alpha)$ can be viewed as an element of $GL(V)$ simply by the trivial extension. In this way one can view $G_0:=\prod_\alpha G_\alpha$ as a Lie subgroup of $G$. Let $Z_0={\rm Gen}\{\exp T_u;\ u\in \bar\frkp_0\}$. Then $Z_0$ is contained in the center $Z(G)$ of $G$ since $\frkp_0$ is contained in the center of $\frkg$. It follows that $Z_0$ can be identified with a commutative Lie group $\{(e^{t_1},\cdots,e^{t_r});\ t=(t_1,\cdots,t_r)\in T_0\}$. Define a subgroup $\bar G$ of $G$ by $\bar G=\{zg;\ z\in Z_0,\ g\in G_0\}\equiv Z_0\times G_0$.

{\lem\label{lem7.1} $\bar G$ is connected, of which the Lie algebra $\cll_{\bar G}$ is identical to $\frkg$. In particular, $G=\bar G$.}

\proof By the definition of $\bar G$ and \eqref{eqn7.2},
$$
\cll_{\bar G}=\cll_{Z_0}\oplus \cll_{G_0} =\frkp_0\oplus (\oplus_\alpha \frkg_\alpha) =\frkp_0\oplus (\oplus_\alpha \frkk_\alpha) \oplus (\oplus_\alpha \frkp_\alpha)=\frkk\oplus \frkp=\frkg.
$$
\endproof

Define $M^n=T_0\times M^{n_1}\times\cdots M^{n_r}$. Then we have

{\prop\label{prop7.2} $M^n=G/K$ as a symmetric space. Furthermore, the equiaffine symmetric hypersurface $x:G/K\to V$ defined by \eqref{eqn7.3} has an expression as follows
\be\label{eqn7.6}\left.\begin{aligned}
x(t,p)=&(c_1e^{t_1}x_1(p_1),\cdots,c_re^{t_r}x_r(p_r)),\\
\forall t=&(t_1,\cdots,t_r)\in T_0,\ p=(p_1,\cdots,p_r)\in M^{n_1}_1\times\cdots\times M^{n_r}_r,\end{aligned}\right.
\ee}
where for each $\alpha=1,\cdots,r$, $c_\alpha$ is a constant uniquely determined by the affine curvatures $\lalp$ and $L_1$ of $x$.

\proof Note that $M^{n_\alpha}_\alpha=G_\alpha/K_\alpha$, $\alpha=1,\cdots,r$ and that $T_0$ is a commutative group with respect to the component-wise addition which is isomorphic to the commutative group $Z_0$ with component-wise multiplication. In this sense,
$$
M^n\cong Z_0\times \prod_\alpha (G_\alpha/K_\alpha)=Z_0\times \left(\left(\prod_\alpha G_\alpha\right)/\left(\prod_\alpha K_\alpha\right)\right)=G/K.
$$

For any $t=(t_1,\cdots,t_r)\in T_0$, $p=(p_1,\cdots,p_r)\in M^n$, we have $p_\alpha=g_\alpha K_\alpha$ for some $g_\alpha\in G_\alpha$, $\alpha=1,\cdots,r$. Note that, as an element of the additive group $T_0$, $t=(t_1,\cdots,t_r)$ is identified with an element $e(t):=(e^{t_1},\cdots,e^{t_r})$ in $Z_0$. Put $g_0=(g_1,\cdots,g_r)\in G_0$. Then $g:=e(t)g_0\in G$ and $(t,p)=gK$.

On the other hand, by the equivariance of $x_\alpha$,
$$x_\alpha(p_\alpha)=x_\alpha(g_\alpha K_\alpha)=g_\alpha x_\alpha(o_\alpha) =C_\alpha g_\alpha\cdot e_\alpha,\quad \alpha=1,\cdots,r.$$
It follows from \eqref{eqn5.7} that
\begin{align*}
x(t,p)=&Cg\cdot e=C(e(t)g_0)\cdot e=C(e^{t_1}g_1\cdot e_1,\cdots,e^{t_r}g_r\cdot e_r)\\
=&C(C^{-1}_1e^{t_1} x_1(p_1),\cdots,C^{-1}_re^{t_r}x_r(p_r))=(c_1e^{t_1}x_1(p_1),\cdots,c_re^{t_r}x_r(p_r)),
\end{align*}
where $c_\alpha=CC^{-1}_\alpha$, $\alpha=1,\cdots,r$, is uniquely determined by $\lalp$ and $L_1$. \endproof

{\dfn\label{dfn7.1}\rm The equiaffine symmetric hypersurface $x:M^n\to V\equiv \bbr^{n+1}$ is called the Calabi composition of the given hypersurfaces $x_\alpha:M^{n_\alpha}_\alpha\to\bbr^{n_\alpha+1}$, $\alpha=1,\cdots,r$.}

{\rmk\rm The Calabi composition formula \eqref{eqn7.6} can be extended to construct new proper affine hyperspheres from more general proper affine hyperspheres including $0$-dimensional ones. Beside the pioneer work in \cite{cal72} by Calabi, other similar composition formulas of Calabi-type can be found in \cite{lix93}, \cite{dil-vra94}, \cite{lix11}, \cite{lix14} and \cite{hil12}.}

\section{The proofs of Theorem \ref{main0}}\label{sec8}

In this section, we shall make use of Theorem \ref{main} and the discussion of the Calabi-type composition in Section \ref{sec7} to provide a proof of Theorem \ref{main0}.

So, let $x:M^n\to \bbr^{n+1}$ be a nondegenerate equiaffine symmetric hypersurface with the affine metric $g$ and affine mean curvature $L_1\neq 0$. Then by Definition \ref{dfn afsym}, $M^n=G/K$ is a symmetric pseudo-Riemannian space with respect to the affine metric $g$, and the Fubini-Pick form $A$ is $G$-invariant. Furthermore, by Proposition \ref{prop3.3}, the symmetric pair $(\frkg,\frkk)$ of Lie algebras can be given by \eqref{eqn3.3} and \eqref{eqn3.4}. Since only the local characterization is considered here, we can assume that $M^n$ is connected and simply connected. It then follows by Theorem \ref{main} that $x$ defines a semi-simple real Jordan algebra ${\mathcal J}=(V,\circ)$ with $V=\bbr^{n+1}$ and, from the structure Lie algebra $\cll$ we have also defined an effective symmetric pair $(\frkg,\frkk)$ by \eqref{eqn4.13} and \eqref{eqn4.23}, where $\frkg$ is the restricted Lie algebra of ${\mathcal J}$.

Now only the following two cases are possible to occur:

(1) The real Jordan algebra ${\mathcal J}=(V,\circ)$ is simple. Note that all the real simple Jordan algebras are listed by the classification theorem (Theorem \ref{thm2.8}). On the other hand, a real split (resp. complex split) quaternion
can be identified with a real (resp. complex) $2\times 2$ matrix (\cite{hil12}, \cite{mcc04})
$$
\lmx a&b\\c&d\rmx,\text{\ with conjugation\ } \ol{\lmx a&b\\c&d\rmx} =\lmx d&-b\\-c&a\rmx,\quad a,b,c,d\in\bbr\text{\ (resp. $\bbc$).}
$$
It follows that the real (resp. complex) Jordan algebra $H_m(Q,\bbr)$ (resp. $H_m(Q,\bbc))$ can be identified with the space of $2m\times 2m$ real (resp. complex) skew-Hamiltonian matrices, that is, matrices in the form
$$
\lmx A&B\\C&A^t\rmx, \text{\ where\ }A,B,C,D\in M_m(\bbr)\ \text{(resp. $M_m(\bbc)$) with\ }B^t=-B,\ C^t=-C
$$
with the Jordan product
$$
X\circ Y=\fr12(XY+YX),\text{\ for all skew-Hamiltonian matrices\ }X,Y
$$
(see \cite{hil12}, \cite{mcc04}). Furthermore, set $J=\lmx 0& I_m\\-I_m&0\rmx$ where $I_m$ is the identity matrix of order $m$. Then, by left multiplication, $J$ isomorphically maps the space of $2m\times 2m$ real (resp. complex) skew-Hamiltonian matrices onto the space $A_{2m}(\bbr)$ (resp. $A_{2m}(\bbc)$) of real (resp. complex) skew-symmetric matrices (\cite{ikr01}); Moreover, if equip $A_{2m}(\bbr)$ (resp. $A_{2m}(\bbc)$) with the Jordan product
$$
X\circ Y=\fr12(XJY+YJX),\quad \forall X,Y\in A_{2m}(\bbr) \text{\ (resp. $A_{2m}(\bbc)$)},
$$
Then $J$ is clearly an isomorphism of Jordan algebras.

Therefore, the list appeared in Theorem \ref{main0} is obtained from Proposition \ref{prop5.3} and the following proposition which can be shown by a long and direct computation case by case, but at the moment we delete this because, in fact, it has already been done by R. Hildebrand (see \cite{hil12}, Section 5):

{\prop (\cite{hil12}) The determinant functions $\det P_u$, $u\in V$, for the simple real Jordan algebras ${\mathcal J}=(V,\circ)$ are as follows:
\begin{enumerate}
\item $V=\bbr$: $\det P_u=u^2$, $u\in \bbr$;
\item $V=Jord_m(Q_\bbr)$, $m\geq 3$: $\det P_u=(u^tQu)^m$, $u\in\bbr^m$;
\item $V=M_m(\bbr)$, $m\geq 3$: $\det P_u=(\det u)^{2m}$, $u\in M_m(\bbr)$;
\item $V=M_m(\bbh)$, $m\geq 2$: $\det P_u=(\det u)^{4m}$, $u\in M_m(\bbh)$;
\item $V=S_m(\bbr,\Gamma)$, $m\geq 3$: $\det P_u=\pm(\det u)^{m+1}$, $u\in S_m(\bbr,\Gamma)$;
\item $V=H_m(\bbc,\Gamma)$, $m\geq 3$: $\det P_u=(\det u)^{2m}$, $u\in H_m(\bbc,\Gamma)$;
\item $V=H_m(\bbh,\Gamma)$, $m\geq 3$: $\det P_u=(\det u)^{2m-1}$, $u\in H_m(\bbh,\Gamma)$;
\item $V=A_{2m}(\bbr)$, $m\geq 3$: $\det P_u=(\det u)^{2m-1}$, $u\in A_{2m}(\bbr)$;
\item $V=SH_m(\bbh)$, $m\geq 2$: $\det P_u=(\det u)^{2m+1}$, $u\in H_m(\bbh,\Gamma)$;
\item $V=H_3(\bbo,\Gamma)$: $\det P_u=(\det u)^{18}$, $u\in H_3(\bbo,\Gamma)$;
\item $V=H_3(O,\bbr)$: $\det P_u=(\det u)^{18}$, $u\in H_3(O,\bbr)$;
\item $V=\bbc$: $\det P_u=|u|^4$, $u\in\bbc$;
\item $V=Jord_m(I)$, $m\geq 3$: $\det P_u=|u^tu|^{2m}$, $u\in\bbc^m$;
\item $V=S_m(\bbc)$, $m\geq 3$: $\det P_u=|\det u|^{2m+1}$, $u\in S_m(\bbc)$;
\item $V=M_m(\bbc)$, $m\geq 3$: $\det P_u=|\det u|^{4m}$, $u\in M_m(\bbc)$;
\item $V=A_{2m}(\bbc)$, $m\geq 3$: $\det P_u=|\det u|^{4m-2}$, $u\in H_m(H,\bbc)$;
\item $V=H_3(O,\bbc)$: $\det P_u=|\det u|^{36}$, $u\in H_3(O,\bbc)$.
\end{enumerate}}

(2) The Jordan algebra ${\mathcal J}$ is not simple. In this case, by Proposition \ref{prop2.2}, ${\mathcal J}$ can be decomposed into a direct sum of some of its simple ideals ${\mathcal J}_\alpha=(V_\alpha,\circ)$ with $n_\alpha=\dim V_\alpha\geq 1$, $\alpha=1,\cdots,r\geq 2$, that is, ${\mathcal J}={\mathcal J}_1\oplus\cdots\oplus {\mathcal J}_r$. By Theorem \ref{main}, each ${\mathcal J}_\alpha$, $1\leq \alpha\leq r$, corresponds one connected, simply connected and nondegenerate hypersurface $x_\alpha:M^{n_\alpha}_\alpha\to\bbr^{n_\alpha+1}$ with affine metric $\galp$, Fubini-Pick form $\Aalp$ and a preassigned affine mean curvature $\lalp\neq 0$ where $n_1+\cdots+n_r=n+1$. Let $\bar x:\bar M^n\to\bbr^{n+1}$ be the Calabi composition of $x_1,\cdots,x_r$ defined in Section \ref{sec7} (see Definition \ref{dfn7.1}). Then by Theorem \ref{main} and the discussion in Section \ref{sec7}, it is not hard to find that $x$ and $\bar x$ are affine equivalent to each other since they correspond to the same real semi-simple Jordan algebra. This completes the proof of Theorem \ref{main0}.\endproof

As the end of the present paper, we give a direct application of the classification theorem (Theorem \ref{main0}) as follows:

Let $x:M^n\to\bbr^{n+1}$ be a nondegenerate hypersurface with affine metric $g$, affine mean curvature $L_1\neq 0$ and a parallel Fubini-Pick form $A$. Then by the affine Gauss equation \eqref{gaus}, the curvature tensor $R$ is parallel with respect to the Levi-Civita connection of $g$, or equivalently, $(M^n,g)$ is a locally pseudo-Riemannian symmetric space. Therefore $(M^n,g)$ is locally isometric to a connected pseudo-Riemannian symmetric space $G/K$ for some symmetric pair $(G,K)$ of groups. On the other hand, by Proposition \ref{prop3.6}, the assumption that $A$ is parallel means that $A$ is $G$-invariant. Thus it follows that $x$ is locally affine equivalent to an equiaffine symmetric hypersurface with nonzero affine mean curvature $L_1$. Therefore, the following theorem is direct by Theorem \ref{main0}:

{\thm\label{thm8.3} (cf. \cite{hil12}) Let $x:M^n\to \bbr^{n+1}$ be a nondegenerate hypersurface with parallel Fubini-Pick form and nonzero affine mean curvature. Then either of the following two cases must occur:

(1) $x$ is locally affine equivalent to a standard imbedding $\td x:G/K\to V$, defined by \eqref{eqn5.7}, of a pseudo-Riemannian symmetric space $G/K$ into an $(n+1)$-dimensional simple Jordan algebra ${\mathcal J}=(V,\circ)$ where the real linear spaces $V$ is among the list in Theorem \ref{main0};

(2) $x$ is locally affine equivalent to a Calabi composition of some of the nondegenerate hypersurfaces in (1) including the $0$-dimensional ones.}

\vskip 0.3in
\flushleft
Xingxiao Li\\
School of Mathematics and Information Sciences\\
Henan Normal University\\
XinXiang 453007, Henan\\
P.R.China\\
email: xxl@henannu.edu.cn

\end{document}